\documentclass{article}%
\usepackage{graphicx}
\usepackage{amsmath}
\usepackage{amsfonts}
\usepackage{amssymb}%
\newtheorem{theorem}{Theorem}

\newtheorem{corollary}[theorem]{Corollary}

\newtheorem{definition}[theorem]{Definition}

\newtheorem{lemma}[theorem]{Lemma}

\newtheorem{proposition}[theorem]{Proposition}

\begin{document}

\title{Configuration of lines in del Pezzo surfaces\\with Gosset Polytopes}
\author{Jae-Hyouk Lee}
\maketitle

\begin{abstract}
In this article, we study the divisor classes$\ $of del Pezzo surfaces, which
are written as the sum of distinct lines with fixed intersection according to
the inscribed simplexes and crosspolytopes in Gosset polytopes.

We introduce the $k$-Steiner system and cornered simplexes, and characterize
the configurations of inscribed $m(<4)$-simplexes with them.

Higher dimensional inscribed $m(4\leq m\leq7)$-simplexes exist in $4_{21}\ $in
the Picard group of del Pezzo surface $S_{8}\ \ $of degree $1.$The
configurations of $4$- and\ $7$-simplexes are related to rulings in $S_{8}$.
And the configurations of $5$- and\ $6$-simplexes correspond the skew
$3$-lines and skew $7$-lines in $S_{8}$. In particular, the seven lines in a
$6$-simplex produce a Fano plane.

We also study the inscribed crosspolytopes and hypercubes in the Gosset
polytopes.\newpage

\end{abstract}
\tableofcontents

\newpage


\section{ Introduction}

In this article, we study the divisor classes$\ $of del Pezzo surfaces which
are written as the sum of distinct lines with fixed intersection, and we
explore the configuration of lines producing the divisors along the geometry
of Gosset polytopes. As the vertices of the Gosset polytopes are corresponded
to lines, the divisor classes we consider are related to polytopes in the
Gosset polytopes. Because Gosset polytopes are compact convex polytopes with
symmetries , divisor classes in this article are naturally finite, and the
corresponding configurations have symmetries leading us interesting results hereinafter.

A \textit{del Pezzo surface}$\ $is a smooth irreducible surface $S_{r}$ whose
anticanonical class $-$ $K_{S_{r}}$ is ample. Each del Pezzo surfaces can be
constructed by blowing up $r\leq$ $8$ points from $\mathbb{P}^{2}$ unless it
is $\mathbb{P}^{1}\times\mathbb{P}^{1}$. A line in $Pic$ $S_{r}\ $is a divisor
class $l\ $with $l^{2}=l\cdot K_{S_{r}}=-1$. By the adjunction formula, a line
$l$ contains a rational smooth curve in $S_{r}$. This rational curve is
embedded to a line in $\mathbb{P}^{9-r}$ along the embedding given by
$\left\vert -K_{S_{r}}\right\vert $, when $S_{r}\ $is very ample. The set of
lines $L_{r}\ $in $Pic\ S_{r}\ $is finite, and its symmetry group is the Weyl
group $E_{r}$.

The \textit{Gosset polytopes} $(r-4)_{21}$, $3\leq r\leq8$, are the
$r$-dimensional semiregular polytopes discovered by Gosset whose symmetry
groups are the Coxeter group $E_{r}$. The vertex figure of $(r-4)_{21}$ is
$\left(  r-5\right)  _{21}$ and the facets of $(r-4)_{21}$ are regular
$(r-1)$-dimensional simplexes $\alpha_{r-1}$ and $(r-1)$-dimensional
crosspolytopes $\beta_{r-1}$. But all the lower dimensional subpolytopes are
regular simplexes. The set of vertices in $(r-4)_{21}\ $is bijective to
$L_{r}\ $the set of lines in $Pic\ S_{r}\ $(\cite{Coxeter4}\cite{Manni}). In
fact, author \cite{Lee} showed that the convex hull of $L_{r}\ $in
$Pic\ S_{r}$ is the Gosset polytopes $(r-4)_{21}$. These are the Gosset
polytopes considered in this article. As the subsets of vertices in
$(r-4)_{21},\ $equivalently lines in $L_{r}$, represent the polytopes in
$(r-4)_{21}$, the divisor classes under our consideration and the
configuration of lines in them are closely related to the the polytopes in n
$(r-4)_{21}$.

The sum of $m$-lines$\ $with $0$-intersection in $Pic\ S_{r}\ $is called a
\textit{skew }$m$\textit{-line} $D_{m}$ in $S_{r}\ $(\cite{Lee}\ or subsection
\ref{Subsec-Gosset-Picard}). According to the correspondence between the lines
in $L_{r}\ $and vertices\ in $(r-4)_{21}$, $L_{r}^{m}\ $the set of skew
$m$-lines in $S_{r}\ $is bijectively related to the set of the regular
$m$-simplexes in $(r-4)_{21}$. Thus the divisor classes $D_{m}$ in
$Pic\ S_{r}\ $exist for $1\leq m\leq r\ $and satisfy $D_{m}^{2}=-m\ $and
$K_{S_{r}}D_{m}=-m$. In fact, skew $m$-lines are equivalent to the divisor
classes with $D^{2}=-m\ $and $K_{S_{r}}D=-m$, when $1\leq m\leq3$. As a line
in $Pic\ S_{r}\ $gives a rational map from $S_{r}\ $to $S_{r-1}\ $by blowing
down the exceptional curve in the line, the skew $r$-lines in $Pic\ S_{r}$ is
corresponded to a rational map from $S_{r}\ $to $\mathbb{P}^{2}$. Therefrom,
we relate the skew $r$-lines in $Pic\ S_{r}\ $to divisor classes $D_{t}\ $with
$D_{t}^{2}=1$ and $K_{S_{r}}D_{t}=-3$. This divisor class in $Pic\ S_{r}\ $is
called an \textit{exceptional system, }and its linear system gives a regular
map to $\mathbb{P}^{2}$. By a transformation given by $K_{S_{r}}+3D_{t}$,
$\mathcal{E}_{r}\ $the set of exceptional system in $Pic\ S_{r}\ $and
$L_{r}^{r}\ $are bijectively related for $3\leq r\leq$ $7$. When $r=8$, the
set of exceptional systems has two orbits of $E_{8}\ $action. One orbit
corresponds to the set of skew $8$-lines in $S_{8}$, and the other orbit
corresponds to the set of $E_{8}$-roots. Thus an exceptional system $D_{t}%
\ $in $\mathcal{E}_{8}\ $satisfies either $-3K_{S_{8}}+2d=D_{t}\ $for a root
$d\ \ $or $3D_{t}+K_{S_{8}}=D_{1}\ $for $D_{1}\ $is a skew $8$-lines. These
two orbits in $\mathcal{E}_{8}\ $play key roles later on.

As the regular $(r-1)$-simplexes, a type of facets in $(r-4)_{21}$,$\ $are
related to exceptional systems, the $(r-1)$-crosspolytopes in $(r-4)_{21}$,
the other type of facets, are also connected to special divisor classes in
$Pic\ S_{r}\ $so called \textit{rulings}. A\ ruling is a divisor class $D$ in
$Pic\ S_{r}$ with $D^{2}=0\ $and $K_{S_{r}}D=-2$ which gives a fibration of
$S_{r}$ over $\mathbb{P}^{1}$. The $F_{r}\ $set of rulings in $Pic\ S_{r}$ is
bijective to the set of $(r-1)$-crosspolytopes in the $(r-4)_{21}%
\ $(\cite{Lee}). Rulings are studied in many different directions. As the
lines play the generators in Cox rings, the rulings determines the relations
for the rings by Batyrev and Popov (\cite{Batyrev-Popov}). And rulings are
also applied to the geometry of the line bundles over del Pezzo surfaces via
the representation theory by Leung and Zhang (\cite{Leung1}\cite{Leung-Zhang}%
).\ Each ruling in $Pic\ S_{r}\ $is bijectively related to the sum of two
lines $L_{r}\ $with $1$-intersection. Here, each of the sum can be written by
$(r-1)\ $different pair of lines which are corresponded to the $(r-1)$-bipolar
pairs of vertices in an $(r-1)$-crosspolytope in $(r-4)_{21}$.

Now, the remain cases of the sums of lines in $L_{r}\ $have positive
intersections. The distinct lines with fixed positive intersections are
corresponded to simplexes whose$\ $vertices are in the vertices of
$(r-4)_{21}\ $but their edge are not. Since $(r-4)_{21}\ $is convex, these
simplexes are inscribed in $(r-4)_{21}$, and the barycenters of the simplexes
are corresponded to the sums of the distinct lines with fixed positive
intersection.\ Therefrom we call the set of $m$-lines $\{l_{1},...,l_{m}%
\}\ $in $L_{r}\ $with fixed intersection $b>0$ as an \textit{inscribed }%
$b$\textit{-degree }$m$\textit{-simplex in }$(r-4)_{21}$ or an $A_{m}^{r}%
(b)$\textit{-polytope}, and the sum of line $l_{1}+...+l_{m}\ $is called as
the center of the $A_{m}^{r}(b)$-polytope. Here the $k$-degree is the common
intersection between vertices, i.e. lines, in the $A_{m}^{r}(b)$-polytope and
also equivalently implies that the length of edges is $\sqrt{2(1+b)}$. As the
maximal possible intersection between lines in $L_{r}\ $is $3$, there is no
inscribed simplexes with degree $m>3$. As the inscribed $1$-degree
$1$-simplexes in $(r-4)_{21}\ $are rulings, we also assume $m>1\ $from now on.
When $r=6$, $A_{2}^{6}(1)$-polytopes are the biggest possible one in $2_{21}$
and there is no $A_{m}^{r}(b)$-polytope with $m>1$.\ When $r=7$, $A_{3}%
^{7}(1)$-polytopes are the biggest possible one in $3_{21}$ with $1$-degree,
and $A_{1}^{7}(2)$-polytopes are the only possible one with $k>2$. When
$r=8,A_{m}^{8}(1)$-polytopes exist up to $m=7$ because the root space of
$S_{8}\ $is $8$-dimensional, and $A_{1}^{8}(2)$,$\ A_{2}^{8}(2)\ $and
$A_{1}^{8}(3)$-polytopes also exist.

For $A_{m}^{r}(1)$-simplexes in $2_{21}$, $3_{21}\ $and $4_{21}\ $for $m\leq
3$, the configurations of lines in the simplexes can be obtained by
$k$-Steiner system and monoidal transforms.

A \textit{Steiner system} $S(x,y,z)\ $is a type of block design system given
by a family of $y$-element subsets, called as blocks, in $z$-element total set
where each $x$-element in the total set is contained exactly one subset in the
family. For example, $S(2,3,7)$ represents the famous Fano projective plane.
The deterministic nature of blocks of Steiner systems can be found in the
Steiner triplet in cubic surfaces which may not be honest Steiner system. But
we focus the deterministic nature of the systems and define $\mathcal{S}%
(k,S_{r})$, $k$\textit{-Steiner system in }$S_{r}$ which is a family of
subsets of $L_{r}\ $where each $(k-1)$-lines in $L_{r}\ $with constant
intersections to each other determines exactly one subset in the family. A
typical type of $k$-Steiner systems is found when the sums of the lines in all
the blocks is constant. In other words, all blocks in the $k$-Steiner system
are corresponded to the $A_{k-1}^{r}(m)$-simplexes for some $m\ $with common
center. In section \ref{Sec-Steiner}, we define $k$-Steiner systems
$\mathcal{S}_{A}(2,S_{7})$,$\ \mathcal{S}_{A}(2,S_{8})$,$\ \mathcal{S}%
_{B}(3,S_{6})$,$\ \mathcal{S}_{B}(3,S_{8})$, and $\mathcal{S}_{C}(4,S_{7}%
)\ $which are related to $A_{1}^{7}(2)$-, $A_{1}^{8}(3)$-, $A_{2}^{6}(1)$-,
$A_{2}^{8}(2)$-\ and $A_{3}^{7}(1)$-polytopes respectively.

To apply the monoidal transform for the inscribed simplexes, we introduce a
notion of cornered simplex. A simplex in $(r-4)_{21}\ $is called
\textit{cornered} if there is a line $l\ $in $L_{r}\ $whose vertex figure
contains the simplex, and the otherwise is called \textit{uncornered}. Here
the \textit{vertex figure of }$l\ $\textit{in }$(r-4)_{21}\ $is the subset of
lines in $L_{r}\ $with $0$-intersection with $l$.\ Thus a simplex cornered by
a line $l\ \ $is preserved by the blow down map $\pi_{l}^{r}:S_{r}\rightarrow
S_{r-1}\ $given by the line $l$. As there is no $A_{m}^{r}(1)$-polytopes in
$2_{21}$, $3_{21}\ $for $m>3$, all the $A_{m>3}^{8}(1)$-polytopes in $4_{21}$
are uncornered. For $A_{3}^{8}(1)$-polytopes in $4_{21}$, both cases appear.
The centers of $A_{3}^{8}(1)$-polytopes in $4_{21}\ $are bijectively related
to the exceptional system in $\mathcal{E}_{8}$, and the centers of cornered
(resp. uncornered) $A_{3}^{8}(1)$-polytopes correspond roots\ (skew $8$-lines)
in $S_{8}$.

The $A_{2}^{r}(1)$-polytopes exist for $r=6,7,8$. All the $A_{2}^{6}%
(1)$-polytopes in $2_{21}\ $share a center so that the configuration of
$A_{2}^{6}(1)$-polytopes in $2_{21}\ $is given by $\mathcal{S}_{B}(3,S_{6}%
)$,$\ 3$-Steiner system in $S_{6}$. The set of $A_{2}^{7}(1)\ $(resp.
$A_{2}^{8}(1)$)-polytopes is equivalently the set of lines in $Pic\ S_{7}%
\ $(resp. skew $2$-lines\ in $Pic\ S_{8}$). The $A_{2}^{7}(1)\ $(resp.
$A_{2}^{8}(1)$)-polytopes are cornered by a line (resp. a skew $2$-line\ in
$Pic\ S_{8}$), and the configuration of $A_{2}^{7}(1)\ $(resp. $A_{2}^{8}%
(1)$)-polytopes\ of a fixed center is also given by $\mathcal{S}_{B}%
(3,S_{6})\ $along the monomial transform.

The $A_{3}^{r}(1)$-polytopes exist for $r=7,8$. All the $A_{3}^{7}%
(1)$-polytopes in $2_{21}\ $share a center so that the configuration of
$A_{3}^{7}(1)$-polytopes in $3_{21}\ $is given by $\mathcal{S}_{C}(4,S_{7}%
)$,$\ 4$-Steiner system in $S_{7}$. The set of cornered $A_{3}^{8}%
(1)$-polytopes is equivalently the set of roots in $Pic\ S_{8}$, and the
configuration of cornered $A_{3}^{8}(1)$-polytopes\ of a fixed center is also
given by $\mathcal{S}_{C}(4,S_{7})\ $along the monomial transform. On the
other hand, the set of uncornered $A_{3}^{8}(1)$-polytopes is equivalent to
the set of $7$-simplexes in $4_{21}$, and the configuration of uncornered
$A_{3}^{8}(1)$-polytopes\ of a fixed center is obtained as the set of $4$-skew
edges in a $7$-simplex in $4_{21}$.

The $A_{m}^{r}(1)$-polytopes for $m>3\ $exist only in $4_{21}$. As
$A_{m>3}^{8}(1)$-polytopes are uncornered, the configuration for these are far
from being uniform. But as the uncornered $A_{3}^{8}(1)$-polytopes\ are
related to the $7$-simplexes in $4_{21}$,the $A_{5}^{8}(1)\ $(resp.$A_{6}%
^{8}(1)$)-polytopes are related to skew $3$-lines (resp. skew $7$-lines), and
the $A_{4}^{8}(1)\ $and $A_{7}^{8}(1)$-polytopes are related to rulings in
$Pic\ S_{8}$.

An $A_{4}^{8}(1)$-polytope\ can be obtained by adding proper line to an
uncornered $A_{3}^{8}(1)$-polytope. And this $A_{4}^{8}(1)$-polytope\ contains
unique corned $A_{3}^{8}(1)$-polytope in it, which gives a line $l$ in
$L_{8}\ $according to the correspondence between the roots (in fact, lines
when $r=8$) in $Pic\ S_{8}\ $and the centers of cornered $A_{3}^{8}%
(1)$\textbf{-}polytope in $4_{21}$. This line $l\ $intersects by one with the
line $l^{\prime}$ remained line in $A_{4}^{8}(1)$-polytope\ after taking off
the corned $A_{3}^{8}(1)$-polytope. Thus\ $l+l^{\prime}\ $gives a ruling in
$Pic\ S_{8}$. It turns out this is true for all the $A_{4}^{8}(1)$-polytopes.
Thus for each center $D_{5}\ $of $A_{4}^{8}(1)$-polytope, there is unique line
$l_{D_{5}}^{S_{8}}\ $in $L_{8}\ $and the center $A_{D_{5}}\ $of an cornered
$A_{3}^{8}(1)$\textbf{-}polytope such that $D_{5}=$ $A_{D_{5}}+l_{D_{5}%
}^{S_{8}}$.\ Furthermore, the set of the centers of $A_{4}^{8}(1)$%
\textbf{-}polytope in $4_{21}\ $is bijective to the following set of order
pairs of lines with $1$-intersection$\ $defined as
\[
\tilde{F}_{8}:=\{(l_{1},l_{2})\mid l_{1}\text{,}l_{2}\in L_{8}\ \text{with
}l_{1}\cdot l_{2}=1\}\text{.}%
\]
Moreover, all the $A_{4}^{8}(1)$\textbf{-}polytopes with center $D_{5}$ in
$4_{21}\ $share a common line $l_{D_{5}}$, and $D_{5}-l_{D_{5}}\ $is the
common center for the unique cornered $A_{3}^{8}(1)$\textbf{-}polytope in
each$\ A_{4}^{8}(1)$\textbf{-}polytope.

The set of centers of $A_{5}^{8}(1)$-polytopes is equivalent to the set of
skew $3$-lines in $Pic\ S_{8}$. As a skew $3$-line is given by unique triple
of lines with $0$-intersection, an $A_{5}^{8}(1)$-polytope contains unique
triple of cornered $A_{3}^{8}(1)$\textbf{-}polytopes in it. This
characterization gives the configuration of $A_{5}^{8}(1)$-polytopes.

Just like $A_{5}^{8}(1)$-polytopes, each center of $A_{6}^{8}(1)$%
\textbf{-}polytope gives a skew $7$-lines in $4_{21}\ $which is given by the
unique choice of seven lines with $0$-intersection, and the uniqueness is the
key to study the configuration of $A_{6}^{8}(1)$\textbf{-}polytopes in
$4_{21}$. Thus for each $A_{6}^{8}(1)$\textbf{-}polytope, there are seven
cornered $A_{3}^{8}(1)$\textbf{-}polytopes in it. Furthermore, for each
cornered $A_{3}^{8}(1)$\textbf{-}polytope in a $A_{6}^{8}(1)$\textbf{-}%
polytope, the remained three lines form an $A_{2}^{8}(1)$\textbf{-}polytope
which is not contained in any$\ $cornered $A_{3}^{8}(1)$\textbf{-}polytope in
the $6$-simplex. We call the triplet of lines \textit{Fano block} and show
that each $A_{6}^{8}(1)$\textbf{-}polytope\ contains seven Fano blocks in it.
Moreover, the seven lines in a $A_{6}^{8}(1)$\textbf{-}polytope\ and its Fano
blocks produce a Steiner system $S(2,3,7)\ $which is known as \textit{Fano
plane}.

The configuration of $A_{7}^{8}(1)$\textbf{-}polytopes in $4_{21}\ $is close
to that of $A_{4}^{8}(1)$\textbf{-}polytopes. Here each center of $A_{7}%
^{8}(1)$\textbf{-}polytopes gives a ruling in $4_{21}$. Furthermore, each
center of $A_{7}^{8}(1)$\textbf{-}polytopes can be written as the sum of two
centers of$\ $cornered $A_{3}^{8}(1)$\textbf{-}polytopes where they are in the
vertex figures of the bipolar pair of lines in the corresponding ruling.

The configurations of inscribed $A_{m}^{r}(b)$\textbf{-}polytopes in
$(r-4)_{21}\ $in degree $b>1\ $are obtained by the $k$-Steiner system and the
monoidal transforms.

We also consider the inscribed crosspolytopes in $(r-4)_{21}\ $which exists
only for $r=8$. And as the applications of the configurations of lines for
$A_{m}^{r}(b)$\textbf{-}polytopes in $(r-4)_{21}\ $and crosspolytopes in
$4_{21}$, we study the hypercubes in $3_{21}\ $and $4_{21}$.

As Coxeter\cite{Coxeter5} related his study on $4_{21}\ $to the Cayley
integral numbers, another direct application of the configuration of lines in
del Pezzo surfaces$\ $can be found in the octonions numbers. We discuss this
in another article.

\section{Preliminaries}

\bigskip

\subsection{\label{Sec-polytopes}Regular Polytopes and Gosset Polytopes}

\bigskip

In this subsection, we review the general theory on regular polytopes that we
use in this article and a family of semiregular polytopes known as Gosset
figures ($k_{21}$ according to Coxeter). Here, we only present general facts,
and the further detail of them can be found \cite{Coxeter}\cite{Coxeter1}%
\cite{Coxeter2}and \cite{Lee}.

Let $P_{n}$\ be a convex $n$-polytope in an $n$-dimensional Euclidean space.
For each vertex $V$, the midpoints of all the edges emanating from a vertex
$V$ in $P_{n}$ form an $(n-1)$-polytope if they lie in a hyperplane, and this
$(n-1)$-polytope is called the \textit{vertex figure} of $P_{n}$ at $V$.\ In
this article, the vertices on the other ends of the edges emanating from the
vertex $V\ $also form an $(n-1)$-polytope, and we also call this
$(n-1)$-polytope the \textit{vertex figure }of $V\ $in $P_{n}$.

A \textit{regular} polytope $P_{n}$ ($n\geq2$) is a polytope whose facets and
vertex figure at each vertex are regular. In particular, a polygon $P_{2}$ is
regular if it is equilateral and equiangular. Naturally, the facets of regular
$P_{n}$ are all congruent, and the vertex figures are all the same.

In this article, we consider two classes of regular polytopes.

\textbf{(1)} A \textbf{regular simplex} $\alpha_{n}$ is an $n$-dimensional
simplex with equilateral edges. Note $\alpha_{n}$ is a pyramid based on
$\alpha_{n-1}$. Thus the facets of a regular simplex $\alpha_{n}$ is a regular
simplex $\alpha_{n-1}$, and the vertex figure of $\alpha_{n}$ is also
$\alpha_{n-1}$. For example, $\alpha_{1}$ is a line-segment, $\alpha_{2}$ is
an equilateral triangle, and $\alpha_{3}$ is a tetrahedron. For a regular
simplex $\alpha_{n}$, only regular simplex $\alpha_{k}$, $0\leq k\leq n-1$
appears as subpolytopes.

\textbf{(2)} A \textbf{crosspolytope} $\beta_{n}$ is an $n$-dimensional
polytope whose $2n$-vertices are the intersects between an $n$-dimensional
Cartesian coordinate frame and a sphere centered at the origin. Note
$\beta_{n}$ is a bipyramid based on $\beta_{n-1}$, and the $n$-vertices in
$\beta_{n}$ form $\alpha_{n-1}$ if the choice is made as one vertex from each
Cartesian coordinate line. So the vertex figure of a crosspolytope $\beta_{n}$
is also a crosspolytope $\beta_{n-1}$, and the facets of $\beta_{n}$ is
$\alpha_{n-1}$. For instance, $\beta_{1}$ is a line-segment, $\beta_{2}$ is a
square, and $\beta_{3}$ is an octahedron. For a crosspolytope $\beta_{n}$,
only regular simplex $\alpha_{k}$, $0\leq k\leq n-1$ appears as subpolytopes.

A polytope $P_{n}$ is called \textit{semiregular} if its facets are regular
and its vertices are equivalent, namely, the symmetry group of $P_{n}$ acts
transitively on the vertices of $P_{n}$.

Here, we consider the semiregular $k_{21}$\ polytopes discovered by Gosset
which are $(k+4)$-dimensional polytopes whose symmetry groups are the Coxeter
group $E_{k+4}$. Note that the vertex figure of $k_{21}$ is $\left(
k-1\right)  _{21}$ and the facets of $k_{21}$ are regular simplexes
$\alpha_{k+3}$ and crosspolytopes $\beta_{k+3}$. The list of $k_{21}$
polytopes is following.%

\[
\underset{\text{List of Gosset Polytopes }k_{21}}{%
\begin{tabular}
[c]{|l|l|l|}\hline
$\ \ k$ & $E_{k+4}$ & $k_{21}$-polytopes\\\hline\hline
$-1$ & $A_{1}\times$$A_{2}$ & triangular prism\\\hline
$\ \ 0$ & $A_{4}$ & rectified 5-cell\\\hline
$\ \ 1$ & $D_{5}$ & demipenteract\\\hline
$\ \ 2$ & $E_{6}$ & $E_{6}$-polytope\\\hline
$\ \ 3$ & $E_{7}$ & $E_{7}$-polytope\\\hline
$\ \ 4$ & $E_{8}$ & $E_{8}$-polytope\\\hline
\end{tabular}
}%
\]

For $k\neq-1$, the facets of $k_{21}$-polytopes are the regular simplex
$\alpha_{k+3}$ and the crosspolytope $\beta_{k+3}$. But all the lower
dimensional subpolytopes are regular simplexes. When $k=-1$, the vertex figure
in $-1_{21}$ is an isosceles triangle instead of an equilateral triangle, and
its facets are the regular triangle $\alpha_{2}$ and the square $\beta_{2}$.

The table of total numbers of subpolytopes in $k_{21}\ $is very useful in this
article and presented below.%

\[
\underset{\ \ \text{Number of subpolytopes in }k_{21}}{%
\begin{tabular}
[c]{|l||l|l|l|l|l|l|}\hline
$E_{k+4}$-polytope($k_{21}$) & $-1_{21}$ & $0_{21}$ & $1_{21}$ & $2_{21}$ &
$3_{21}$ & $4_{21}$\\\hline\hline
$\beta_{k+3}$ & $\ \ 3$ & $5$ & $10$ & $27$ & $126$ & $2160$\\\hline
vertex & $\ \ 6$ & $10$ & $16$ & $27$ & $56$ & $240$\\\hline
$\alpha_{1}$ & $\ \ 9$ & $30$ & $80$ & $216$ & $756$ & $6720$\\\hline
$\alpha_{2}$ & $\ \ 2$ & $30$ & $160$ & $720$ & $4032$ & $60480$\\\hline
$\alpha_{3}$ &  & $5$ & $120$ & $1080$ & $10080$ & $241920$\\\hline
$\alpha_{4}$ &  &  & $16$ & $648$ & $12096$ & $483840$\\\hline
$\alpha_{5}$ &  &  &  & $72$ & $6048$ & $483840$\\\hline
$\alpha_{6}$ &  &  &  &  & $576$ & $207360$\\\hline
$\alpha_{7}$ &  &  &  &  &  & $17280$\\\hline
\end{tabular}
\ }%
\]

\bigskip

\subsection{\label{Subsec-Gosset-Picard}Gosset Polytopes in the Picard Groups
of del Pezzo Surfaces}

\bigskip

The del Pezzo surfaces are smooth irreducible surfaces $S_{r}$ whose
anticanonical class $-$ $K_{S_{r}}$ is ample. We can construct the del Pezzo
surfaces by blowing up $r\leq$ $8$ points from $\mathbb{P}^{2}$ unless it is
$\mathbb{P}^{1}\times\mathbb{P}^{1}$. In particular, it is very well known
that there are $27$ lines on a cubic surface $S_{6}$ and the configuration of
these lines is acted by the Weyl group $E_{6}$(\cite{Demaz}\cite{Dolgachev}%
\cite{Hartsh}). The set of $27$-lines in $S_{6}$ are bijective to the set of
vertices of a Gosset $2_{21}$ polytope. The similar correspondences are found
between the $28$-bitangents in $S_{7}\ $and $3_{21}$ polytopes, and between
the tritangent planes for $S_{8}$ and $4_{21}$ polytopes. The correspondence
between lines in $S_{6}$ and vertices in $2_{21}$ is applied to study the
geometry of $2_{21}$ by Coxeter(\cite{Coxeter4}), and the correspondence is
extended to each $3\leq r\leq$ $8$(\cite{Manni}).

We denote such a del Pezzo surface by $S_{r}$ and the corresponding blow up by
$\pi_{r}:$ $S_{r}\rightarrow$ $\mathbb{P}^{2}$. And $K_{S_{r}}^{2}$ $=9-r$ is
called the degree of the del Pezzo surface. Each exceptional curve and the
corresponding class given by blowing up is denoted by $e_{i}$, and both the
class of $\pi_{r}^{\ast}\left(  h\right)  $ in $S_{r}$ and the class of a line
$h$ in $\mathbb{P}^{2}$ are referred as $h$. Then, we have
\[
h^{2}=1\text{, }h\cdot e_{i}=0\text{, }e_{i}\cdot e_{j}=-\delta_{ij}%
\ \ \text{for }1\leq i,j\leq r,
\]
and the Picard group\ of $S_{r}$ is
\[
Pic\ S_{r}\simeq\mathbb{Z}h\oplus\mathbb{Z}e_{1}\oplus...\oplus\mathbb{Z}e_{r}%
\]
with the signature $(1,-r)$. And $K_{S_{r}}=-3h+\sum_{i=1}^{r}e_{i}$.

The inner product given by the intersection on $Pic\ S_{r}\ $induces a
negative definite metric on $\left(  \mathbb{Z}K_{S_{r}}\right)  ^{\perp}$ in
$Pic\ S_{r}\ $where we can also define natural reflections. To define
reflections on $\left(  \mathbb{Z}K_{S_{r}}\right)  ^{\perp}$ in $Pic\ S_{r}$,
we consider a root system%

\[
R_{r}:=\{d\in Pic\ S_{r}\mid d^{2}=-2,\ d\cdot K_{S_{r}}=0\}\text{,}%
\]
with simple roots $d_{0}=h-e_{1}-e_{2}-e_{3},d_{i}=e_{i}-e_{i+1},\ 1\leq i\leq
r-1.\ $Each element $d$ in $R_{r}$ defines a reflection on $\left(
\mathbb{Z}K_{S_{r}}\right)  ^{\perp}$ in $Pic\ S_{r}$
\[
\sigma_{d}(D):=D+\left(  D\cdot d\right)  d\ \ \text{for\ }D\in\left(
\mathbb{Z}K_{S_{r}}\right)  ^{\perp}%
\]
and the corresponding Weyl group $W(S_{r})$ is $E_{r}$ where $3\leq r\leq8$.
Furthermore, the reflection $\sigma_{d}$ on $\left(  \mathbb{Z}K_{S_{r}%
}\right)  ^{\perp}\ $can be used to obtain a transformation both on
$Pic\ S_{r}$ and on $Pic\ S_{r}\otimes\mathbb{Q\simeq Q}h\oplus\mathbb{Q}%
e_{1}\oplus...\oplus\mathbb{Q}e_{r}$ via the linear extension of the
intersections of divisors in $Pic\ S_{r}$. Here $Pic\ S_{r}\otimes\mathbb{Q}$
is a vector space with the signature $(1,-r)$.

In this article, we deal with divisor classes $D\ $satisfying $D\cdot
K_{S_{r}}=\alpha$,\ $D^{2}=\beta\ $where $\alpha\ $and $\beta\ $are integers,
which are preserved by the extended action of $W(S_{r})$. In particular,
$W(S_{r})\ $acts as a reflection group on the set of divisor classes with
$D\cdot K_{S_{r}}=\alpha$. Therefrom, we define an \textit{affine hyperplane
section} in $Pic\ S_{r}\otimes\mathbb{Q}$ defined by
\[
\tilde{H}_{b}:=\{D\in Pic\ S_{r}\ \otimes\mathbb{Q}\mid-D\cdot K_{S_{r}}=b\}
\]
where $b$ is an arbitrary real number and an affine hyperplane section
$H_{b}:=\tilde{H}_{b}\cap Pic\ S_{r}\ $in $Pic\ S_{r}$. By the ample condition
of $-K_{S_{r}}\ $and the Hodge index theorem, the inner product on
$Pic\ S_{r}$ induces a negative definite metric on $H_{b}$. In fact, the
induced metric is defined on $Pic\ S_{r}\otimes\mathbb{Q}$, and we can also
consider the induced norm by fixing a center $\frac{b}{9-r}K_{S_{r}}$ in the
affine hyperplane section $-D\cdot K_{S_{r}}=b$ in $Pic\ S_{r}\otimes
\mathbb{Q}$. This norm is also negative definite. Furthermore, the reflection
$\sigma_{d}\ $on $Pic\ S_{r}\otimes\mathbb{Q\ }$induces a reflection on
$\tilde{H}_{b}$.\ Generically, the hyperplanes in $Pic\ S_{r}\otimes
\mathbb{Q\ }$induce affine hyperplanes in $\tilde{H}_{b}\ $and they may not
share a common point. But the reflection hyperplane of each reflection
$\sigma_{d}\ $in $Pic\ S_{r}\otimes\mathbb{Q\ }$gives a hyperplane in
$\tilde{H}_{b}\ $containing the center because it is given by a condition
$K_{S_{r}}\cdot d=0$. Thus, the Weyl group $W(S_{r})\ $acts on $\tilde{H}%
_{b}\ $and $H_{b}\ $as an reflection group.

Now, we want construct Gosset polytopes $(r-4)_{21}$ in $Pic\ S_{r}%
\otimes\mathbb{Q}\ $as the convex hull of the set of special classes in
$Pic\ S_{r}$, which is known as lines. A $line\ $in $Pic\ S_{r}\ $is
equivalently a divisor class $l\ $with $l^{2}=-1\ $and $K_{S_{r}}\cdot l=-1$,
and the set of lines is given as
\[
L_{r}:=\{D\in Pic(S_{r})\mid D^{2}=-1,K_{S_{r}}D=-1\}\text{.}%
\]
As the Weyl group $W(S_{r})\ $acts as an affine reflection group on the affine
hyperplane given by $D\cdot K_{S_{r}}=-1$, $W(S_{r})\ $acts on the set of
lines in $Pic\ S_{r}$. Therefrom, we construct a semiregular polytope in
$Pic\ S_{r}\otimes\mathbb{Q}\ $whose vertices are exactly the lines in
$Pic\ S_{r}$.\ Since the symmetry group of the polytope is $W(S_{r})$, the
polytope is actually a Gosset polytope $(r-4)_{21}$.

\bigskip

\textbf{Remark: }Each line $l$ in $L_{r}\ $contains a exceptional curve which
produces a blow down map from $S_{r}\ $to $S_{r-1}$. We denote the map as
$\pi_{l}^{r}:S_{r}\rightarrow S_{r-1}$.

\bigskip

For a Gosset polytope $(r-4)_{21}$, subpolytopes are regular simplexes excepts
the facets which consist of$\ (r-3)$-simplexes and $(r-3)$-crosspolytopes.
Since the subpolytopes in $(r-4)_{21}$ are basically configurations of
vertices, we obtain natural characterization of subpolytopes in $(r-4)_{21}%
\ $as divisor classes in $Pic\ S_{r}$.

\bigskip

\textbf{Remark : }To identify each subpolytope in $(r-4)_{21}$, we want to use
the barycenter of the subpolytope. Each vertex of the polytope $(r-4)_{21}$
represents a line in $S_{r}$, and the honest centers of simplexes (resp.
crosspolytopes ) are written as $(l_{1}+...+l_{k})/k$ (resp.$(l_{1}^{\prime
}+l_{2}^{\prime})/2$)\ in $\tilde{H}_{1}$ which may not be elements in $Pic$
$S_{r}$. Therefore, alternatively, we choose $(l_{1}+...+l_{k})$ as the center
of a subpolytope so that $(l_{1}+...+l_{k})$ is in $Pic$ $S_{r}$.

\bigskip

We use the algebraic geometry of del Pezzo surfaces to identify the divisor
classes corresponding to the subpolytopes in $(r-4)_{21}$. For this purpose,
we consider the following set of divisor classes which are called\textit{
}$skew$\textit{ }$a$\textit{-}$lines$\textit{, }$exceptional$\textit{
}$systems$ and\ \textit{rulings }in $Pic\ S_{r}$.%

\begin{align*}
L_{r}^{a}  &  :=\{D\in Pic(S_{r})\mid D=l_{1}+...+l_{a},\text{ }%
l_{i}\ \text{disjoint lines in }S_{r}\}\\
\mathcal{E}_{r}  &  :=\{D\in Pic(S_{r})\mid D^{2}=1,\text{ }K_{S_{r}}\cdot
D=-3\}\\
F_{r}  &  :=\{D\in Pic(S_{r})\mid D^{2}=0,\text{ }K_{S_{r}}\cdot D=-2\}
\end{align*}

The \textit{skew }$a$\textit{-line} in $L_{r}^{a}$ is an extension of the
definition of lines in $S_{r}$. Each skew $a$-line represents an
$(a-1)$-simplex in an $(r-4)_{21}$ polytope. Furthermore, each skew $a$-line,
there is only one set of disjoint lines in $L_{8}$ to present it. The skew
$a$-lines also holds $D^{2}=-a$ and $D\cdot K_{S_{r}}=-a$, and the divisor
classes with these conditions are equivalently skew $a$-lines for $a\leq3$.

The \textit{exceptional system} in $\mathcal{E}_{r}\ $is a divisor class in
$Pic\ S_{r}$ whose linear system gives a regular map from $S_{r}$ to
$\mathbb{P}^{2}$. As this regular map corresponds to a blowing up from
$\mathbb{P}^{2}$ to $S_{r}$, naturally exceptional systems are related to
$(r-1)$-simplexes in $(r-4)_{21}$ polytopes, which is the one of two types of
facets appearing in $(r-4)_{21}$ polytopes. In fact, by a transformation
$\Phi\ $from $\mathcal{E}_{r}\ $to $L_{r}^{r}\ $by
\[
\Phi(D_{t}):=K_{S_{r}}+3D_{t}\ \text{for }D_{t}\in\mathcal{E}_{r}\text{,}%
\]
the set $\mathcal{E}_{r}$ is bijective to the set of the $(r-1)$-simplexes in
$(r-4)_{21}$ polytopes, for $3\leq r\leq$ $7$. When $r=8$, the set of
exceptional systems has two orbits. One orbit with $17280$ elements
corresponds to the set of skew $8$-lines in $S_{8}$, and the other orbit with
$240$ elements corresponds to the set of $E_{8}$-roots because $-3K_{S_{8}}+$
$2d$ is an exceptional system for each $E_{8}$-root $d$. Thus an exceptional
system $D_{t}\ $in $\mathcal{E}_{r}\ $satisfies either $-3K_{S_{8}}%
+2d=D_{t}\ $for a root $d\ \ $or $3D_{t}+K_{S_{8}}=D_{1}\ $for $D_{1}\ $is a
skew $8$-lines.

The \textit{ruling }in $F_{r}$ is a divisor class in $Pic\ S_{r}$ which gives
a fibration of $S_{r}$ over $\mathbb{P}^{1}$. And we show that the $F_{r}$ is
bijective to the set of $(r-1)$-crosspolytopes in the $(r-4)_{21}$ polytope.
Furthermore, we explain the relationships between lines and rulings according
to the incidence between the vertices and $(r-1)$-crosspolytopes. This leads
us that a pair of proper crosspolytopes in the $(r-4)_{21}$ give the blowing
down maps from $S_{r}$ to $\mathbb{P}^{1}\times$ $\mathbb{P}^{1}$.

After proper comparison between divisor classes obtained from the geometry of
the polytope $(r-4)_{21}\ $and those given by the geometry of a del Pezzo
surface, we come by the following correspondences.%

\[%
\begin{tabular}
[c]{|l|l|}\hline
del Pezzo surface $S_{r}$ & Gosset polytopes $(r-4)_{21}$\\\hline\hline
lines & vertices\\\hline
skew $a$-lines $1\leq a\leq r$ & $\left(  a-1\right)  $-simplexes $1\leq a\leq
r$\\\hline
exceptional systems & $\left(  r-1\right)  $-simplexes ($r<8$)\\\hline
rulings & $\left(  r-1\right)  $-crosspolytopes\\\hline
\end{tabular}
\ \
\]

\textbf{Lines and rulings}

Now, we know rulings in $Pic\ S_{r}\ $correspond $\left(  r-1\right)
$-crosspolytopes in the Gosset polytopes $(r-4)_{21}$. Since an $\left(
r-1\right)  $-crosspolytope\ has $2(r-1)\ $vertices in it, we want to have a
criterion that a line $l\ $in $L_{r}\ $is one of the vertices of an $\left(
r-1\right)  $-crosspolytope\ corresponding to a ruling $f$. Here, the
following are equivalently stating the criterion.

(1) $f\cdot l=0$

(2) $f-l$ is a line

(3) The vertex$\ $represented by $l$ in $(r-4)_{21}$\ is one of the vertex of
the $(r-1)$-crosspolytope corresponding to $f$.

\bigskip

\textbf{Gieser Transform and Bertini Transform for lines}

In \cite{Lee}, we define the \textit{Gieser transform on lines in }$S_{7}$ or
simply Gieser transform by%
\[
G(l):=-\left(  K_{S_{7}}+l\right)  \ \text{for }l\in L_{7}\text{.}%
\]
This is an automorphism of $L_{7}\ $which is closely related to the deck
transformation of the covering from $S_{7}$ to $\mathbb{P}^{2}$ given by
$\left\vert -K_{S_{7}}\right\vert $. Similarly a transformation $B$ on lines
in $S_{8}$ defined by
\[
B(l):=-\left(  2K_{S_{8}}+l\right)  \ \text{for\ }l\in L_{8}%
\]
is referred as the \textit{Bertini transform on lines} or simply Bertini transform.

The Gieser transform $G$ and the Bertini transform $B$ preserve the
intersection between two lines. Therefrom, these transforms can be extended to
symmetries on the configurations of lines in $S_{7}\ $and $S_{8}$,\ which are
also symmetries on $3_{21}$ and $4_{21}$ respectively (see \cite{Lee}).

\subsection{\textbf{\label{Subsec-monoidal}}Monoidal transform of del Pezzo
surfaces}

Since a del Pezzo surface $S_{r}$ is obtained by blowing up of one point on
$S_{r-1}$, we can describe divisor classes in$\ S_{r-1}\ $producing lines in
$S_{r}$ after blowing up. Let $l$ be the exceptional divisor in $S_{r}$ given
by a blow-up of a point from $S_{r-1}$, namely $\pi_{l}^{r}:S_{r}\rightarrow
S_{r-1}$. The proper transform of a divisor $D$ in $Pic\ S_{r-1}$ producing a
line in $S_{r}\ $satisfies
\[
(\pi_{l}^{r\ast}\left(  D\right)  -ml)^{2}=-1\text{, }(\pi_{l}^{r\ast}\left(
D\right)  -ml)\cdot(\pi_{l}^{r\ast}\left(  K_{S_{r-1}}\right)  +l)=-1
\]
for a nonnegative integer $m$. Therefore, we consider a divisor $D$ in
$Pic(S_{r-1})$ with
\[
D^{2}=m^{2}-1\text{, }D\cdot K_{S_{r-1}}=-m-1.
\]
By the Hodge index theorem, the list of possible $m$ is
\[
m=\left\{
\begin{tabular}
[c]{l}%
$0,1\ $\ $\ \ \ \ \ \ \ \ $\\
$0,1,2\ \ \ \ $\\
$0,1,2,3\ $%
\end{tabular}%
\begin{tabular}
[c]{l}%
$4\leq r\leq6$\\
$r=7$\\
$r=8$%
\end{tabular}
\ \ \ \ \right.  .
\]

On the other hand, the integer $m\ $is the intersection between $l\ $and a
line $\pi_{l}^{r\ast}\left(  D\right)  -ml$. Therefore, the above divisors in
$S_{r-1}\ $characterize the subsets of $L_{r}\ $according to the intersection
with $l$. Here we consider a set $N_{k}(l,S_{r})$ defined as%
\[
N_{k}(l,S_{r}):=\left\{  \ l^{\prime}\ \in L_{r}\mid l^{\prime}\cdot
l=m\right\}  \ \text{for }k\geq-1\text{,}%
\]
and obtain the following bijections.%
\begin{align*}
N_{0}(l,S_{r})  &  \approx L_{r-1}\ \text{for }4\leq r\leq8\\
N_{1}(l,S_{r})  &  \approx F_{r-1}\ \text{for }4\leq r\leq8\\
N_{2}(l,S_{8})  &  \approx N_{0}(-2K_{S_{8}}-l,S_{8})\approx L_{7}%
\text{,\ }N_{2}(l,S_{7})=\{-K_{S_{7}}-l\}\text{\ }\\
N_{3}(l,S_{8})  &  =\{-2K_{S_{8}}-l\}\text{.}%
\end{align*}

It will be useful to observe the Gieser transform $G\ $on $L_{7}\ $gives
\[
G(N_{-1}(l,S_{7}))=N_{2}(l,S_{7})\text{ and\ }G(N_{0}(l,S_{7}))=N_{1}%
(l,S_{7})\text{,}%
\]
and the Bertini transform $B\ $on $L_{8}\ $induces%
\[
B(N_{-1}(l,S_{8}))=N_{3}(l,S_{8})\text{, }B(N_{0}(l,S_{8}))=N_{2}%
(l,S_{8})\ \text{and }B(N_{1}(l,S_{8}))=N_{1}(l,S_{8})\text{.}%
\]

\bigskip

\textbf{Lines of vertex figure }

In this article, the subset $N_{0}(l,S_{r})\ $in $L_{r}\ $plays very important
role. When we consider a rational map given by blowing down of the exceptional
curve in $l$, $N_{0}(l,S_{r})\ $in $L_{r}\ $is bijectively mapped to the lines
in $S_{r-1}\ $which are also vertices of the corresponding Gosset polytope
$(r-5)_{21}$. We call $N_{0}(l,S_{r})\ $in $L_{r}$%
\ $\mathit{lines\ of\ vetrex\ figure\ }$of $l$ or simply vertex figure of
$l$.\ In fact, each line in $N_{0}(l,S_{r})\ $is joined to $l\ $with an edge,
and the set of midpoints of there edges is the vertex figure of the vertex
$l\ $in $(r-4)_{21}$.

\bigskip

\section{\label{Sec-Steiner}Steiner Systems of del Pezzo Surfaces}

\bigskip

A \textit{Steiner system} $S(a,b,c)$ is a type of block design system which is
a family of subsets $S\ $consisting of $b$-elements in a set $T$ of
$c$-elements satisfying that each $a$-elements of $T$ is contained exactly one
subset in the family. Here the subset $S\ $in the system is called a
\textit{block}.

For example, $S(2,3,7)$ represents the famous Fano projective plane which is a
family of $3$-points subsets (called lines) in $7$-points set$\ $where each
two points in the set determines a line. This is one special case of Steiner
triple systems $S(2,3,n)$ containing $n(n-1)/6$ blocks. An analog of Steiner
system can be consider on the set of lines in the Picard group of $S_{6}$
whose block is given as a subset(called a triplet) consisting of three lines
in $S_{6}$ with $1$-intersection to each others. Because the sum of these
three lines in a block equals to $-K_{S_{6}}$, a block is defined by each two
lines with $1$-intersection. However, the total number of blocks in this
system is much less than expected number of blocks in the ordinary Steiner
triple system $S(2,3,27)$. Therefore, it is natural to add a few combinatorial
conditions to define an analogue of Steiner system on the family of divisors
of del Pezzo surfaces. In this section, we define Steiner systems on del Pezzo
surfaces as an analogous to the triplets in cubic surfaces and search Steiner
systems in the configurations of lines in del Pezzo surfaces.

In this article, we work on divisor classes $D\ $given as a sum of lines with
fixed intercessions. The divisor class $D\ $satisfies equations $D^{2}%
=\alpha\ $and $D\cdot K_{S_{8}}=\beta$. For certain integers $\alpha\ $and
$\beta$, the equations have unique solution which gives a condition
determining a block in the following $k$-Steiner system.

\begin{definition}
A family of subset of lines in a del Pezzo surface $S_{r}$ written
$\mathcal{S}(k,S_{r})\ \ $is called a $k$-Steiner system on $S_{r}\ $if it is
a family of subset $S\ $of $k$-lines in $Pic$ $S_{r}\ $where each
$(k-1)$-lines in $Pic$ $S_{r}\ $with constant intersection to each other
determines exactly a subset $S\ $in $\mathcal{S}(k,S_{r})$.
\end{definition}

\bigskip Note: Each block in $\mathcal{S}(k,S_{r})\ $has $k$-lines.

\bigskip

\textbf{A. }$2$\textbf{-Steiner system on }$S_{7}\ $\textbf{and} $S_{8}$

We consider a family of subsets of lines in $S_{7}\ $defined by
\[
\mathcal{S}_{A}(2,S_{7}):=\{\{l_{1},l_{2}\}\mid l_{1}\cdot l_{2}%
=2\ \ \text{for\ }l_{1}\text{, }l_{2}\ \in L_{7}\}\text{.}%
\]
Recall that for two lines $l_{1}$ and $l_{2}$ in $Pic\ S_{7}$, $l_{1}\cdot
l_{2}=2$ is equivalent to $l_{1}+l_{2}=-K_{S_{7}}$. Thus each line $l$ in
$S_{7}$ gives a $2$-Steiner block defined as $\{l,-\left(  K_{S_{7}}+l\right)
\}=\{l,G\left(  l\right)  \}$. Therefore, $\mathcal{S}_{A}(2,S_{7})\ $is a
$2$-Steiner system on $S_{7}$. In fact, each $2$-Steiner block in
$\mathcal{S}_{A}(2,S_{7})\ $on $S_{7}$ represents a bitangent of degree $2$
covering from $S_{7}$ to $\mathbb{P}^{2}$ given by $\left\vert -K_{S_{7}%
}\right\vert $ (see chapter 8 \cite{Dolgachev}).

Similarly, we define a family of subsets of lines in $S_{8}\ $by%
\[
\mathcal{S}_{A}(2,S_{8}):=\{\{l_{1},l_{2}\}\mid l_{1}\cdot l_{2}%
=3\ \ \text{for\ }l_{1}\text{, }l_{2}\ \in L_{8}\}\text{.}%
\]
Because any two lines $l_{1}$ and $l_{2}$ in $Pic\ S_{8}$ with $l_{1}\cdot
l_{2}=3$ equivalently hold $l_{1}+l_{2}=-2K_{S_{8}}$.\ Therefore, each line
$l$ in $S_{8}$ determine a $2$-Steiner block given by $\{l,-\left(  2K_{S_{8}%
}+l\right)  \}=\{l,B\left(  l\right)  \}$. So $\mathcal{S}_{A}(2,S_{8})\ $is a
$2$-Steiner system on $S_{8}$. Again, this $2$-Steiner block on $S_{8}$ is
related to a tritangent plane of degree $2$ covering from $S_{8}$ to
$\mathbb{P}^{2}$ given by $\left\vert -2K_{S_{8}}\right\vert $ (see chapter 8
\cite{Dolgachev}).

\bigskip

\textbf{B. }$3$\textbf{-Steiner system on }$S_{6}\ $\textbf{and} $S_{8}$

The triplet in cubic surfaces is a typical example of a $3$-Steiner system on
$S_{6}\ $defined by%
\[
\mathcal{S}_{B}(3,S_{6}):=\{\{l_{1},l_{2},l_{3}\}\mid l_{1}\cdot l_{2}%
=l_{2}\cdot l_{3}=l_{3}\cdot l_{1}=1\ \ \text{for\ }l_{1}\text{, }%
l_{2}\text{,}l_{3}\in L_{6}\}\text{,}%
\]
because each pair of lines $l_{1}\ $and $l_{2}\ $with $l_{1}\cdot l_{2}%
=1\ $determine a line $l:=-K_{S_{6}}-l_{1}-l_{2}\ $with $l_{1}\cdot
l=l_{2}\cdot l=1$.

Similarly, for $S_{8}$%
\[
\mathcal{S}_{B}(3,S_{8}):=\{\{l_{1},l_{2},l_{3}\}\mid l_{1}\cdot l_{2}%
=l_{2}\cdot l_{3}=l_{3}\cdot l_{1}=2\ \ \text{for\ }l_{1}\text{, }%
l_{2}\text{,}l_{3}\in L_{8}\}
\]
is a $3$-Steiner system on $S_{8}\ $because each three line $l_{1}$, $l_{2}$
and $l_{3}$ in $S_{8}$ with $l_{i}\cdot l_{j}=2$, $i\not =j$, holds
\[
l_{1}+l_{2}+l_{3}=-3K_{S_{8}}.
\]
Here, $l_{1}+l_{2}+l_{3}$ is a divisor $D$ with $D^{2}=9$ and $D\cdot
K_{S_{8}}=-3$ and moreover $-3K_{S_{8}}$ is the only divisor with the
conditions\ by Hodge index theorem.

\bigskip

\textbf{Remark}: For each line $l$ in $S_{8}$, $l+K_{S_{8}}$ is a root in
$S_{8}$, and the vice versa.\ Thus,\ any three roots $d_{1}$, $d_{2}$ and
$d_{3}$ in $S_{8}$ with $d_{i}\cdot d_{j}=1$, $i\not =j$ have
\[
d_{1}+d_{2}+d_{3}=l_{1}+l_{2}+l_{3}+3K_{S_{8}}=0,
\]
and also give a $3$-Steiner system of the roots on $S_{8}$.

\bigskip

As an application of this $3$-Steiner system on $S_{8}$, we have the following
theorem. We also provide another proof given by monoidal transform of lines.
The argument in the second proof is very useful in this article.

\begin{theorem}
\label{Thm-K8}For a del Pezzo surface $S_{8},$ the canonical class $K_{S_{8}}$
can be written as
\[
K_{S_{8}}=l_{1}-l_{2}-l_{3}%
\]
where $l_{1},l_{2}\ $and $l_{3}\ $are lines with $l_{1}\cdot l_{2}=l_{1}$
$\cdot l_{3}=0$ and $l_{2}$ $\cdot l_{3}=2$.
\end{theorem}

\textbf{Proof(1)}: Each lines $l_{2}$ and $l_{3}$ with $l_{2}$ $\cdot l_{3}=2$
determines $\ $a$\ $line $-3K_{S_{8}}-l_{2}-l_{3}$ as in a $3$-Steiner block.
Again, as in the $2$-Steiner system, we get a line
\[
l_{1}=-\left(  2K_{S_{8}}+\left(  -3K_{S_{8}}-l_{2}-l_{3}\right)  \right)
=K_{S_{8}}+(l_{2}+l_{3}).
\]
This gives theorem.

\textbf{Proof(2)}: Consider a line $l_{1}\ $in $L_{8}\ $and corresponding blow
down map $\pi_{l_{1}}^{8}:S_{8}\rightarrow S_{7}$.\ Recall that $-K_{S_{7}}%
\ $can be written as the sum of two lines $l_{2}\ $and $l_{3}\ $in $L_{7}%
\ $with $l_{2}\cdot l_{3}=2$. Since $K_{S_{8}}=\pi_{l_{1}}^{8\ast}(K_{S_{7}%
})+l_{1}\ $and $l_{2}\cdot l_{1}=l_{3}\cdot l_{1}=0$, we have
\[
K_{S_{8}}=-l_{2}-l_{3}+l_{1}\ \text{.}%
\]

${\LARGE \blacksquare}$

\textbf{Remark : }The theorem implies that for a divisor class $l_{2}+l_{3}%
\ $with $l_{2}$ $\cdot l_{3}=2$,\ the divisor class $K_{S_{8}}+l_{2}+l_{3}%
\ $is\ a line with $0$-intersection to $l_{2}\ $and $l_{3}$.

\bigskip

\textbf{C. }$4$\textbf{-Steiner system on }$S_{7}$

For $S_{7}$, any three lines in $S_{7}$ with $1$-intersection to each other
determine a line because each four lines $l_{1}$, $l_{2}$, $l_{3}$ and $l_{4}$
in $S_{7}$ with $l_{i}\cdot l_{j}=1$, $i\not =j$, holds
\[
l_{1}+l_{2}+l_{3}+l_{4}=-2K_{S_{7}}.
\]
Here, $l_{1}+l_{2}+l_{3}+l_{4}$ is a divisor $D$ with $D^{2}=8$ and $D\cdot
K_{S_{7}}=-4$ and $-2K_{S_{7}}$ is the only divisor with the conditions.

Therefore, we get a $4$-Steiner system in $S_{7}\ $given by%
\[
\mathcal{S}_{C}(4,S_{7}):=\{\{l_{1},l_{2},l_{3},l_{4}\}\mid l_{i}\cdot
l_{j}=1,i\not =j\ \ \text{for\ }l_{1}\text{, }l_{2}\text{,}l_{3}\text{,}%
l_{4}\in L_{7}\}\text{.}%
\]

\bigskip

\section{Inscribed regular polytopes in Gosset polytopes}

\bigskip

In this section, we study the convex regular polytopes whose vertices are the
subset of vertices of Gosset polytopes. As the nature of the study on Gosset
polytopes including the symmetry of high degree, the characterization of
inscribed polytopes is complicated. Since we are dealing the divisor classes
written as sum of lines with fixed intersection, we focus on the fundamental
regular polytopes such as simplexes and crosspolytopes in this article. In
fact, one can also consider polytopes, such as cubes, sharing not only
vertices but also edges with Gosset polytopes. We will discuss these in the
next section.

As in \cite{Lee}, we identify the barycentric centers of simplexes and
crosspolytopes inscribed in the Gosset polytopes and study the configuration
of those polytopes. Since we are working along the configuration of lines in
del Pezzo surface, the centers are chosen in $Pic\ S_{r}\ $instead of
$Pic\ S_{r}\otimes\mathbb{Q}$ after multiplying proper integers.

\subsection{Inscribed $1$-degree simplexes\ in Gosset polytopes}

\bigskip

We observe that for any two distinct lines $l_{1}$ and $l_{2}$ in $S_{r}$, we
have
\[
\left(  l_{1}-l_{2}\right)  ^{2}=-2-2l_{1}\cdot l_{2},\ \left(  l_{1}%
-l_{2}\right)  \ \cdot K_{S_{r}}=0.
\]
And since the metric on $K_{S_{r}}^{\perp}$ is negative definite, $l_{1}%
-l_{2}$ has length$\sqrt{2\left(  1+l_{1}\cdot l_{2}\right)  }$. The line
segment joining lines $l_{1}$ and $l_{2}\ $which are two vertices in a Gosset
polytope in $Pic\ S_{r}$ are called an $m$\textit{-degree edge} if $l_{1}\cdot
l_{2}=m$. The $0$-degree edges are honest edges appearing in the Gosset
polytopes. Since the Gosset polytopes are convex, any polytope whose vertices
is the subset of the vertices of the Gosset polytopes is inscribed in them.
Here we consider the following inscribed polytopes.

\begin{definition}
For a Gosset polytope $(r-4)_{21}$ of a del Pezzo surface $S_{r}$, an
$n$\textit{-simplex whose vertex set is the subset of the vertex set of
}$(r-4)_{21}\ $is called an inscribed $m$\textit{-degree }$n$\textit{-simplex}
in $(r-4)_{21}$ if each edge of the simplex has $m$-degree, and we denote it
by $A_{n}^{r}(m)$-polytope.
\end{definition}

\bigskip

As the monoidal transform is useful in this article, we want to separate
inscribed polytopes in $(r-4)_{21}\ $which are preserved by a blow down from
$S_{r}\ $to $S_{r-1}$. Thus this type of inscribed polytopes can be treated
not only in $(r-4)_{21}\ $but also in $(r-5)_{21}$.

\begin{definition}
An inscribed simplex in a Gosset polytope $(r-4)_{21}$ of a del Pezzo surface
$S_{r}$ is called \textit{cornered} if there is a line $l\ $in $L_{r}\ $where
all the lines representing vertices of the simplex$\ $are in the vertex figure
of $l$. Otherwise, the simplex is called uncornered.
\end{definition}

\bigskip

\subsubsection{Inscribed $1$-degree $1$-simplexes and $2$-simplexes}

\bigskip

The centers of inscribed $m(\leq3)$-simplexes in $(r-4)_{21}\ $are either
unique or corresponded to the lines for the cornered issues. Thus the
configurations of inscribed $m(\leq3)$-simplexes in $(r-4)_{21}\ $are
naturally related to the $k$-Steiner systems along the monoidal transform.

\textbf{A. Inscribed }$1$\textbf{-degree }$1$\textbf{-simplexes}

The center of an$\ A_{1}^{r}(1)$-polytope in $(r-4)_{21}$ for $3\leq r\leq
8\ $is the divisor class $D\ $with $D\cdot K_{S_{r}}=-2\ $and $D^{2}=0$, which
is a ruling. Thus by subsection \ref{Subsec-Gosset-Picard}, the set of the
centers of $A_{1}^{r}(1)$-polytopes in $(r-4)_{21}$ is bijective to $F_{r}$
the set of rulings in $S_{r}$. Furthermore, the configuration of the
$A_{1}^{r}(1)$-polytopes with a fixed center $D\ \ $is the set of antipodal
pairs of lines in a crosspolytope in $(r-4)_{21}\ $whose center is $D$.

\bigskip

\textbf{B. Inscribed }$1$\textbf{-degree }$2$\textbf{-simplex}

The $1$-degree inscribed $2$-simplexes exist when $6\leq r\leq8$.

The center of each $A_{2}^{r}(1)$-polytope represents a divisor class
$D:=l_{1}+l_{2}+l_{3}$ where $l_{1},l_{2}$ and $l_{3}$ are lines with
$l_{i}\cdot l_{j}=1,$ $i\not =j$, which satisfies $D^{2}=\left(  l_{1}%
+l_{2}+l_{3}\right)  ^{2}=3$ and $D\cdot K_{S_{r}}=\left(  l_{1}+l_{2}%
+l_{3}\right)  \cdot K_{S_{r}}=-3$.

\bigskip

\textbf{(1) Inscribed }$1$\textbf{-degree }$2$\textbf{-simplex\ in }$2_{21}$

The class \textbf{ }$-K_{S_{6}}$ is the only divisor class satisfying the
above equations for the center. Thus, all the $A_{2}^{6}(1)$-polytopes in
$2_{21}$ share one center. In fact, the triplets in cubic surfaces given by
$1$-intersecting $3$-lines correspond to $A_{2}^{6}(1)$-polytopes in $2_{21}$.
Thus the configuration of $A_{2}^{6}(1)$-polytope in $2_{21}\ $is given by the
$3$-Steiner system$\ \mathcal{S}_{B}(3,S_{6})\ $in section \ref{Sec-Steiner}.

\bigskip

\textbf{(2) Inscribed }$1$\textbf{-degree }$2$\textbf{-simplex in }$3_{21}$

For each divisor class $D\ $satisfying the above equations for the center, we
have a new divisor given by $l_{D}^{S_{7}}:=D+K_{S_{7}}$. This divisor
$l_{D}^{S_{7}}\ $is a line in $L_{7}\ $because $\left(  l_{D}^{S_{7}}\right)
^{2}=$ $\left(  D+K_{S_{7}}\right)  ^{2}=-1\ $and $l_{D}^{S_{7}}\cdot
K_{S_{r}}=\left(  D+K_{S_{7}}\right)  \cdot K_{S_{r}}=-1$. Therefore, the set
of lines in $S_{7}$ is bijective to the set of centers of $A_{2}^{7}%
(1)$-polytopes in $3_{21}$.

Suppose $l_{1}$, $l_{2}\ $and $l_{3}\ $are three vertices of an $A_{2}^{7}%
(1)$-polytope$\ $with the center $D$. Then the line $l_{D}^{S_{7}}=D+K_{S_{7}%
}\ $has
\[
l_{D}\cdot l_{i}=(l_{1}+l_{2}+l_{3}+K_{S_{7}})\cdot l_{i}=0\ \ \text{for
}i=1,2,3\text{.}%
\]
Thus $l_{1}$, $l_{2}\ $and $l_{3}\ $are in $N_{0}(l_{D},S_{7})$, namely,
vertex figure of $l_{D}^{S_{7}}$, and the simplex $A_{2}^{7}(1)\ $is cornered.
Furthermore, $l_{1}$, $l_{2}\ $and $l_{3}\ $induce vertices of an $A_{2}%
^{6}(1)$-polytope in $2_{21}\ $by the blow down map $\pi_{l_{D}^{S_{7}}}%
^{7}:S_{7}\rightarrow S_{6}$. Therefore, the configuration of $A_{2}^{7}%
(1)$-polytopes with fixed center in $3_{21}$ equals the $3$-Steiner
system$\ \mathcal{S}_{B}(3,S_{6})\ $in $S_{6}$.

\bigskip

\textbf{Remark :} The choice of line $l_{D}^{S_{7}}\ $is natural once we
observe its role in the following configuration of lines in $3_{21}$.\ Recall
each two pair of $l_{1}$, $l_{2}\ $and $l_{3}\ $in an $A_{2}^{7}(1)$-polytope
forms a ruling corresponding to a crosspolytope (see subsection
\ref{Subsec-Gosset-Picard} or \cite{Lee}). As the crosspolytopes in $3_{21}%
\ $has $6$-pair of bipolar vertices in it, a ruling can be obtained from
$6$-pair of lines with intersection $1$. We choose $l_{1}+l_{2}\ $for a
ruling, then $G(l_{3})\ $and $l_{D}^{S_{7}}\ $are another bipolar pair of
lines in the ruling because $G(l_{3})+l_{D}^{S_{7}}=\left(  -K_{S_{7}}%
-l_{3}\right)  +\left(  l_{1}+l_{2}+l_{3}+K_{S_{7}}\right)  =l_{1}+l_{2}$.
From the other two choices of rulings from $l_{1}$, $l_{2}\ $and $l_{3}\ $in
the $A_{2}^{7}(1)$-polytope,\ we have similar results containing $l_{D}%
^{S_{7}}\ $in common.\ The line $l_{D}^{S_{7}}\ $is the common vertex of the
crosspolytopes corresponding to the three rulings. Thus $l_{1}$, $l_{2}\ $and
$l_{3}\ $are in the vertex figure of $l_{D}^{S_{7}}$.

\bigskip

\textbf{(3) Inscribed }$1$\textbf{-degree }$2$\textbf{-simplex in }$4_{21}$

Let $D\ $be a divisor class$\ $representing a center of an $A_{2}^{8}%
(1)$-polytope in $4_{21}$, and we consider a divisor defined by $D+K_{S_{8}}$.
Then $D+K_{S_{8}}\ $is a skew $2$-line in $Pic\ S_{8}\ $because $\left(
D+K_{S_{8}}\right)  ^{2}=-2\ $and $\left(  D+K_{S_{8}}\right)  \cdot K_{S_{8}%
}=-2$. Therefore, the set of skew $2$-lines in $S_{8}$ is bijective to the set
of centers of an $A_{2}^{8}(1)$-polytope in $4_{21}$.

Assume $l_{1}$, $l_{2}\ $and $l_{3}\ $form an $A_{2}^{8}(1)$-polytope in
$4_{21}\ $where its center is $D_{1}=l_{1}+l_{2}+l_{3}$.$\ $Then$\ D_{1}%
+K_{S_{8}}\ $is a skew $2$-line which can be written as a sum of two disjoint
lines, namely, $D_{1}+K_{S_{8}}=l_{a}^{123}+l_{b}^{123}\ $where $l_{a}%
^{123}\ $and $l_{b}^{123}\ $in $L_{8}\ $with $l_{a}^{123}\ \cdot l_{b}%
^{123}=0$. Here the choice of $l_{a}^{123}\ $and $l_{b}^{123}\ $is unique
since a skew $2$-line is given by a unique pair of lines. Furthermore we have
\[
\left(  l_{a}^{123}+l_{b}^{123}\right)  \cdot l_{i}=\left(  l_{1}+l_{2}%
+l_{3}+K_{S_{8}}\right)  \cdot l_{i}=0\ \text{for }i=1,2,3\text{.}%
\]

If $l_{a}^{123}\cdot l_{1}=-1\ $and $l_{b}^{123}\cdot l_{1}=1$, then
$l_{a}^{123}=l_{1}\ $and $1=l_{b}^{123}\cdot l_{1}=l_{b}^{123}\cdot
l_{a}^{123}=0$ which is a contradiction. Thus we obtain that $l_{a}^{123}\cdot
l_{i}=l_{b}^{123}\cdot l_{i}=0\ $for $i=1,2,3$, and moreover $l_{1}$,
$l_{2}\ $and $l_{3}\ $are in the vertex figures of $l_{a}^{123}\ $and
$l_{b}^{123}$. Here the simplex is cornered. Since the skew $2$-line
$l_{a}^{123}+l_{b}^{123}\ $induce a blow down map from $S_{8}\ $to $S_{6}$,
and it sends the $A_{2}^{8}(1)$-polytope in $4_{21}$ given by $l_{1}$,
$l_{2}\ $and $l_{3}\ $to an $A_{2}^{6}(1)$-polytope in $2_{21}$. Thus the
configuration of the $A_{2}^{8}(1)$-polytopes in $4_{21}\ $with fixed center
equals to the $3$-Steiner system$\ \mathcal{S}_{B}(3,S_{6})\ $in $S_{6}$.

\bigskip

The configuration of inscribe $1$-degree $2$-simplexes are summarized as follows.

\begin{theorem}
The $3$-Steiner system$\ \mathcal{S}_{B}(3,S_{6})\ $in $S_{6}\ $determines the
configuration of $A_{2}^{8}(1)$-polytopes in $4_{21}$\ with a fixed center,
the configuration of $A_{2}^{7}(1)$--polytopes$\ $in $3_{21}$ with fixed
center and the configuration of $A_{2}^{6}(1)$-polytopes$\ $in $2_{21}$.
\end{theorem}

\bigskip

\subsubsection{Inscribed $1$-degree $3$-simplexes}

\bigskip

The $1$-degree inscribed $3$-simplexes, $A_{3}^{r}(1)$\textbf{-}polytopes,
exist when $r=7,8$.\ The center of each $A_{3}^{r}(1)$\textbf{-}polytope
in\ $(r-4)_{21}\ $represents a divisor class $D:=$ $l_{1}+l_{2}+l_{3}+l_{4}$
where $l_{1},l_{2},l_{3}$ and $l_{4}$ lines with $l_{i}\cdot l_{j}=1$
$i\not =j$, which satisfies $D^{2}=\left(  l_{1}+l_{2}+l_{3}+l_{4}\right)
^{2}=8$ and $D\cdot K_{S_{r}}=$ $\left(  l_{1}+l_{2}+l_{3}+l_{4}\right)  \cdot
K_{S_{r}}=-4$.

\bigskip

\textbf{A. Inscribed }$1$\textbf{-degree }$3$\textbf{-simplex in }$3_{21}$

By Hodge index theorem, $-2K_{S_{7}}$ is the only divisor class satisfying
equations for the center of $A_{3}^{7}(1)$\textbf{-}polytopes in\ $3_{21}$.
Thus all $A_{3}^{7}(1)$\textbf{-}polytopes for $r=7\ $share a common center
$-2K_{S_{7}}$. Furthermore, since $A_{3}^{7}(1)$\textbf{-}polytopes in
$3_{21}$ corresponds to $4$-Steiner blocks in $S_{7}$, the configuration of
$A_{3}^{7}(1)$\textbf{-}polytopes in $3_{21}$ equals the $4$-Steiner
$\mathcal{S}_{C}(4,S_{7})\ $on $S_{7}$ in the section \ref{Sec-Steiner}.

\bigskip

\textbf{B. Inscribed }$1$\textbf{-degree }$3$\textbf{-simplex\ in }$3_{21}$

Each divisor class $D\ $in $Pic\ S_{8}\ $satisfying above equations for center
can be can be transformed to an exceptional system $D+K_{S_{8}}$ because
$\left(  D+K_{S_{8}}\right)  ^{2}=1$ and $\left(  D+K_{S_{8}}\right)  \cdot
K_{S_{8}}=-3$. According to the subsection \ref{Subsec-Gosset-Picard},
$D+K_{S_{8}}$ can be in one of two orbits which correspond to the set of roots
and the set of skew $8$-lines in $Pic\ S_{8}$. If $\left(  l_{1}+l_{2}%
+l_{3}+l_{4}\right)  +K_{S_{8}}$ is corresponded to a root, there is a root
$d$ such that $\left(  l_{1}+l_{2}+l_{3}+l_{4}\right)  +K_{S_{8}}=-3K_{S_{8}%
}+2d\ \ $from the subsection \ref{Subsec-Gosset-Picard}. (1)We consider
\begin{align*}
\left(  l_{1}+l_{2}+l_{3}+l_{4}\right)   &  =\left(  (h-e_{1}-e_{2}%
)+(h-e_{3}-e_{4})+(h-e_{5}-e_{6})+(h-e_{7}-e_{8})\right) \\
&  =h-K_{S_{8}}.
\end{align*}
And $(h-K_{S_{8}})+4K_{S_{8}}=$ $h+3K_{S_{8}}$ cannot be $2d$ for any root
$d$. $\ $(2) However if we choose four lines
\begin{align*}
\left(  l_{1}+l_{2}+l_{3}+l_{4}\right)   &  =\left(
\begin{array}
[c]{c}%
e_{1}+(h-e_{1}-e_{2})+(2h-e_{1}-e_{3}-e_{4}-e_{5}-e_{6})\\
+(3h-e_{1}-e_{2}-e_{3}-e_{4}-e_{5}-e_{6}-2e_{7})
\end{array}
\right) \\
&  =-2K_{S_{8}}+2e_{8},
\end{align*}
then we have $(-2K_{S_{8}}+2e_{8})+4K_{S_{8}}=$ $2(K_{S_{8}}+e_{8})$ where
$(K_{S_{8}}+e_{8})$ is a root.

Therefore, the center $\left(  l_{1}+l_{2}+l_{3}+l_{4}\right)  $ can be
corresponded to an element in the set of skew $8$-lines in $S_{8}$ or the set
of roots in $S_{8}$. Furthermore, by $E_{8}$ action, \ the set of centers of
$A_{3}^{8}(1)$-polytopes in $4_{21}$ is bijective to the set of exceptional
systems on $S_{8}$ which is the union of the two sets which are bijectively
corresponded to the set of skew $8$-lines and the set of roots in $S_{8}$.

\bigskip

\textbf{(1) Cornered inscribed }$1$\textbf{-degree }$3$\textbf{-simplex in
}$4_{21}$

Before we describe the configurations of $A_{3}^{8}(1)$-polytopes in $4_{21}$,
we need to find a better way to separate the$\ A_{3}^{8}(1)$\textbf{-}%
-polytopes in $4_{21}\ $into the above two different sets. Here we consider
the condition \textit{cornered}.\ If an $A_{3}^{8}(1)$\textbf{-}polytope$\ $in
$4_{21}\ $is cornered, there is a line $l\ $where all the lines in an
$A_{3}^{r}(1)$\textbf{-}polytope$\ $are in the vertex figure of $l$. Since
lines of $S_{8}$ correspond to roots of $S_{8}$,\ the line $l\ $whose vertex
figure containing the cornered $A_{3}^{8}(1)$\textbf{-}polytope$\ $gives a
root $l+K_{S_{8}}$. And this roughly implies that cornered $A_{3}^{8}%
(1)$\textbf{-}polytopes in $4_{21}\ $correspond to the roots of $S_{8}$. In
fact, we can explain explicitly the correspondence in the following proposition.

\bigskip

\begin{proposition}
\label{Prop-root-3simplex}Let $D\ $be a center of an $A_{3}^{8}(1)$%
\textbf{-}polytope$\ $in $4_{21}$. The divisor class $D\ $corresponds to a
root $d$ of $S_{8}\ $as $D+4K_{S}=2d\ $if and only if the $A_{3}^{8}%
(1)$\textbf{-}polytope$\ $is cornered.
\end{proposition}

\textbf{Proof}: Suppose $D\ $gives a root $d\ $by $D+4K_{S}=2d$. Then we can
define a line $l:=d-K_{S_{8}}=D/2+K_{S_{8}}$.$\ $For each line $l^{\prime}%
\ $of the $A_{3}^{8}(1)$\textbf{-}polytope$\ $in $4_{21}\ $with center $D$, we
have $l^{\prime}\cdot l=l^{\prime}\cdot\left(  D/2+K_{S_{8}}\right)  =0$.
Therefore, all the lines consisting of the $A_{3}^{8}(1)$\textbf{-}%
polytope$\ $are in the vertex figure of $l$, and the $A_{3}^{8}(1)$%
\textbf{-}polytope$\ $is cornered.

Assume the $A_{3}^{8}(1)$\textbf{-}polytope$\ $is cornered, then there is a
line $l\ $where each lines in the $A_{3}^{8}(1)$\textbf{-}polytope$\ $is in
the vertex figure of $l$. Furthermore, the divisor class $D\ $and the line
$l\ $are related by $D=2l-2K_{S_{8}}$. This relation can be verified by
checking $(2l-2K_{S_{8}}-D)^{2}=0\ $and $(2l-2K_{S_{8}}-D)\cdot K_{S_{8}}%
=0\ $with $D\cdot l=0$. Now the divisor $d\ $given by $l\ +K_{S_{8}}$
satisfies
\[
2d=2(l+K_{S_{8}})=\left(  2l-2K_{S_{8}}\right)  +4K_{S_{8}}=D+4K_{S}\text{.}%
\]

${\LARGE \blacksquare}$

\bigskip

\textbf{Remark :} The cornered condition in the proposition is determined by
the center of the simplex rather than the set of vertices. Therefore, all the
$A_{3}^{8}(1)$\textbf{-}-polytopes in $4_{21}\ $with the center $D\ $%
corresponds to a root $d$ of $S_{r}\ $are cornered.

\bigskip

\begin{corollary}
\label{CoroCornered-line}For a cornered $A_{3}^{8}(1)$\textbf{-}polytope in
$4_{21}\ $given by lines $l_{1},l_{2},l_{3}\ $and $l_{4}$, the line $l\ $in
the proposition is written as
\[
l=K_{S_{8}}+\frac{\left(  l_{1}+l_{2}+l_{3}+l_{4}\right)  }{2},
\]
and $l\ $is the only line in $L_{8}\ $where its vertex figure contains
$l_{1},l_{2},l_{3}\ $and $l_{4}$.
\end{corollary}

\textbf{Proof: }The expression of $l\ $is from the proposition
\ref{Prop-root-3simplex}, and we only need to prove the uniqueness of line
$l$. Assume $l_{a}\ $is another line in $L_{8}\ $where $\{l_{1},l_{2}%
,l_{3},l_{4}\}\subset N_{0}(l_{a},S_{8})$. By checking the possible
intersection $l\cdot l_{a}$, we show that $l=l_{a}$.

The intersection $l_{a}\cdot l\ $is not zero because two skew lines $l_{a}$
and $l\ $produce a blow down from $S_{8}\ $to $S_{6}\ $which induces an
inscribed $1$-degree $3$-simplex in $2_{21}$.\ If $l_{a}\cdot l\ =1$,
$l_{a}+l\ $is a ruling on $S_{8}$ where the lines $l_{1},l_{2},l_{3}\ $and
$l_{4}\ $are the vertices of the $7$-crosspolytope in $4_{21}\ $corresponding
to $l_{a}+l\ $(see subsection \ref{Subsec-Gosset-Picard}). But a line in a
$7$-crosspolytope has only one other line in the $7$-crosspolytope with
$1$-intersection.\ Thus $l_{a}\cdot l\ \not =1$. If $l_{a}\cdot l\ =2,\ $by
theorem \ref{Thm-K8}$\ \ l_{a}+l+K_{S_{8}}\ $gives another line $l_{b}\ $such
that $l_{i}\cdot l_{b}=l_{i}\cdot\left(  l_{a}+l+K_{S_{8}}\right)  =-1\ $for
$i=1,2,3,4$. Thus $l_{i}=l_{b}\ $for $i=1,2,3,4\ $and this is a contradiction.
Therefore $l_{a}\cdot l\ \not =2$. If $l_{a}\cdot l\ =3$, then $l_{a}%
=B(l)=-2K_{S_{8}}-l$. But since $l_{1}\cdot l=0$, we have a contradiction that
$l_{a}\cdot l_{1}=\left(  -2K_{S_{8}}-l\right)  \cdot l_{1}=2$. At last,
$l_{a}\cdot l\ \not =3$. By subsection \ref{Subsec-monoidal}, the intersection
of two lines in $L_{8}\ $is no bigger than $3$. Thus from above, $l_{a}\cdot
l\ $must be $-1$. Thus $l=l_{a}$.

${\LARGE \blacksquare}$

\bigskip

As in subsection \ref{Subsec-Gosset-Picard}, because the set of the centers of
the cornered $A_{3}^{8}(1)$\textbf{-}polytopes$\ $in $4_{21}\ $leads us to the
set of roots in $S_{8}$, we expect the set of the centers of the uncornered
$A_{3}^{8}(1)$\textbf{-}polytopes$\ $in $4_{21}\ $corresponds to the set of
skew $8$-lines in $Pic\ S_{8}$. The proof of the following corollary is
similar to the proposition \ref{Prop-root-3simplex}, and we leave it to the readers.

\begin{corollary}
\label{Coro-Skew8-cornered}Let $D\ $be a center of an $A_{3}^{8}(1)$%
\textbf{-}polytope$\ $in $4_{21}$. The divisor class $D\ $corresponds to a
skew $8$-lines $D_{1}$ in $L_{1}^{8}\ $as $3D+4K_{S_{8}}=D_{1}\ $if and only
if the $A_{3}^{8}(1)$\textbf{-}polytope is uncornered.
\end{corollary}

\textbf{Remark :} This corollary implies that all the $A_{3}^{8}(1)$%
\textbf{-}polytopes in $4_{21}\ $with the center $D\ $corresponds to a skew
$8$-line, namely a $7$-simplex,$\ $are uncornered.

\bigskip

From above corollary \ref{CoroCornered-line}, we can obtain the following
important proposition.

\begin{proposition}
\label{Prop-4cornerTo3+1Uncornered}Let $l_{i}$,$\ 1\leq i\leq5$, be lines in
$L_{8}\ $such that $l_{i}\cdot l_{j}=1\ $for $i\not =j$. If $l_{1},l_{2}%
,l_{3}\ $and $l_{4}\ $form a cornered $A_{3}^{8}(1)$\textbf{-}polytope in
$4_{21}$, another $A_{3}^{8}(1)$\textbf{-}polytope\ given by $l_{1}%
,l_{2},l_{3}\ $and $l_{5}\ $is uncornered.
\end{proposition}

\textbf{Proof: }Let $l\ $be the line where $l_{1},l_{2},l_{3}\ $and $l_{4}%
\ $are in the vertex figure of $l$. By corollary \ref{CoroCornered-line},
$l\ $is $K_{S_{8}}+\frac{l_{1}+l_{2}+l_{3}+l_{4}}{2}$. Now we assume the
$A_{3}^{r}(1)$\textbf{-}polytope\ given by $l_{1},l_{2},l_{3}\ $and $l_{5}%
\ $is cornered. Then there is a line $l^{\prime}\ $given by $K_{S_{8}}%
+\frac{l_{1}+l_{2}+l_{3}+l_{5}}{2}$. We consider
\[
l\cdot l^{\prime}=\left(  K_{S_{8}}+\frac{l_{1}+l_{2}+l_{3}+l_{4}}{2}\right)
\cdot\left(  K_{S_{8}}+\frac{l_{1}+l_{2}+l_{3}+l_{5}}{2}\right)  =-\frac{1}%
{2}.
\]
Since the intersection must be an integer, this is a contradiction to the
assumption. This gives the proposition.\ 

${\LARGE \blacksquare}$

\bigskip

Now the configuration of cornered $A_{3}^{8}(1)$\textbf{-}polytopes$\ $in
$4_{21}\ $is as follows.

Let $D\ $be a center of a cornered $A_{3}^{8}(1)$\textbf{-}polytope$\ $in
$4_{21}$. Then we get a line $l=D/2+K_{S_{8}}\ $where the simplex is in the
vertex figure of $l$. By a blow down map $\pi_{l}^{8}:S_{8}\rightarrow S_{7}$
given by $l$, the simplex in $4_{21}\ $induce an $A_{3}^{7}(1)$\textbf{-}%
polytope$\ $in $3_{21}$, which is a $4$-Steiner\ block in $\mathcal{S}%
_{C}(4,S_{7})\ $on $S_{7}$. Therefore, the configuration of cornered
$A_{3}^{8}(1)$\textbf{-}polytopes$\ $in $4_{21}\ $with a fixed center equals
the $4$-Steiner\ system in $\mathcal{S}_{C}(4,S_{7})\ $on $S_{7}$.

\bigskip

By combining this and the result in from the case \textbf{ }$r=7$\textbf{ }of
$A_{3}^{7}(1)$\textbf{-}polytopes, we have the following theorem.

\begin{theorem}
Both the configuration of cornered $A_{3}^{8}(1)$\textbf{-}polytopes$\ $in
$4_{21}\ $with a fixed center and the configuration of $A_{3}^{7}%
(1)$\textbf{-}polytopes$\ $in $3_{21}$ equal the $4$-Steiner\ system in
$\mathcal{S}_{C}(4,S_{7})\ $on $S_{7}$
\end{theorem}

\bigskip

\textbf{(2) Uncornered inscribed }$1$\textbf{-degree }$3$\textbf{-simplex in
}$4_{21}$

Recall that for each center$\ $of $A_{2}^{8}(1)$\textbf{-}polytopes in
$4_{21}$gives a skew $2$-line which corresponds an edge in $4_{21}$. Since a
skew $2$-line is given by a unique pair of lines, the center also determine
two disjoint lines. By using this fact, we have the following facts for the
configuration of uncornered $A_{3}^{8}(1)$\textbf{-}polytopes in $4_{21}$.

\bigskip

\begin{definition}
Two edges $D_{1}\ $and $D_{2}$, namely, skew $2$-lines, in $(r-4)_{21}\ $of a
del Pezzo surface $S_{r}\ $are called \textit{skew edges} if they are two
disjoint edges in a ($0$-degree) $2$-simplex, equivalent, two disjoint edges
with $D_{1}\cdot$ $D_{2}=0$.
\end{definition}

\bigskip

\begin{proposition}
\label{Prop-skewedges}Let $D$ be a center of an $A_{3}^{8}(1)$\textbf{-}%
polytope in $4_{21}$. If the inscribed simplex whose center $D\ $is
uncornered, $4$-lines in the inscribed simplex induce $4$-skew edges in a
$7$-simplex in $4_{21}\ $whose center is $3D+4K_{S_{8}}$.
\end{proposition}

\textbf{Proof: }Let $l_{1}$, $l_{2}$, $l_{3}\ $and $l_{4}\ $be four lines in
the $A_{3}^{8}(1)$\textbf{-}polytope in $4_{21}\ $with center $D=l_{1}%
+l_{2}+l_{3}+l_{4}$. There are four $A_{2}^{8}(1)$\textbf{-}polytopes in
$4_{21}\ $given by each three choice from $\{l_{1},l_{2},l_{3},l_{4}\}$. For
example, we consider an $A_{2}^{8}(1)$\textbf{-}polytope in $4_{21}\ $given by
$l_{1}$, $l_{2}\ $and $l_{3}$ with the center $l_{1}+l_{2}+l_{3}$. Then there
is a skew $2$-lines\ $l_{a}^{123}+l_{b}^{123}\ $given by $l_{1}+l_{2}%
+l_{3}+K_{S_{8}}=l_{a}^{123}+l_{b}^{123}$ , and $l_{1}$, $l_{2}\ $and
$l_{3}\ $are in the vertex figures of $l_{a}^{123}\ $and$\ l_{a}^{123}$. But
$l_{i}^{123}\cdot l_{4}\geq1\ $for $i=a,b\ $because the $3$-simplex is
uncornered. Furthermore, we have $l_{i}^{123}\cdot l_{4}=1\ $for
$i=a,b\ \ $since%
\[
2=\left(  l_{1}+l_{2}+l_{3}+K_{S_{8}}\right)  \cdot l_{4}=\left(  l_{a}%
^{123}+l_{b}^{123}\right)  \cdot l_{4}=l_{a}^{123}\cdot l_{4}+l_{b}^{123}\cdot
l_{4}\text{.}%
\]
Thus from the lines $l_{1}$, $l_{2}\ $and $l_{3}$ which form an$\ A_{2}%
^{8}(1)$\textbf{-}polytope in the uncornered $A_{3}^{8}(1)$\textbf{-}polytope
in $4_{21}$, we have two disjoint lines $l_{a}^{123}\ $and$\ l_{b}^{123}%
\ $such that%
\[
l_{j}\cdot l_{i}^{123}=0\ \text{and }l_{4}\cdot l_{i}^{123}=1\ \text{for
}i=1,2,3\ \text{and }i=a,b\text{.}%
\]

By performing this process, we have four pair of edges $l_{a}^{lmn}%
+l_{b}^{lmn}\ $for $lmn\in S=:\{123,124,134,234\}$ in $4_{21}$. Here all the
eight lines $l_{i}^{lmn}$ must be distinct because of their intersections to
$l_{1}$, $l_{2}$, $l_{3}\ $and $l_{4}$. Therefore, we have four disjoint edges
in $4_{21}$. Furthermore, we observe
\[
\left(  l_{a}^{234}+l_{b}^{234}\right)  \cdot\left(  l_{a}^{134}+l_{b}%
^{134}\right)  =\left(  l_{2}+l_{3}+l_{4}+K_{S_{8}}\right)  \cdot\left(
l_{1}+l_{3}+l_{4}+K_{S_{8}}\right)  =0
\]
and conclude the four disjoint edges are skew edges to each other. This
implies the distinct eight lines $l_{i}^{lmn}$\ have $0$-intersection to each
other.\ Therefore, they are the vertices of a $7$-simplex in $4_{21}$.
Moreover, the center $D_{1}\ $of the $7$-simplex is
\[
D_{1}=\sum_{lmn\ \in S}(l_{a}^{lmn}+l_{b}^{lmn})=3\sum_{j=1}^{4}%
l_{j}+4K_{S_{8}}\text{.}%
\]

This completes the proposition.$\ \ \ \ \ \ \ \ \ {\LARGE \blacksquare}$

\textbf{Remark:} This proposition shows how to identify the skew $8$-lines in
the corollary \ref{Coro-Skew8-cornered}.

As the converse of above proposition, we have the following theorem.

\begin{theorem}
\label{Thm-uncorner-skewedge}Let $D_{1}\ $be a center of a $7$-simplex in
$4_{21}$. For each $4$-skew edges in a $7$-simplex in $4_{21}$, there is an
$A_{3}^{8}(1)$\textbf{-}polytope in $4_{21}$,$\ $which is uncornered.
Moreover, each edge in the set give a line in the inscribed simplex.
\end{theorem}

\textbf{Proof: }Let $a_{1}$, $a_{2}$, $a_{3}\ $and $a_{4}\ $are $4$-skew edges
in a $7$-simplex in $4_{21}\ $whose center is $a_{1}+a_{2}+a_{3}+a_{4}=D_{1}$.
We denote $a_{i}=l_{i}^{1}+l_{i}^{2}\ $for $i=1,2,3,4\ $where $l_{i}^{1}%
$,$\ l_{i}^{2}\ $in $L_{8}$.\ Since we have skew edges, $l_{i}^{j}$'s\ are
distinct $8$-lines which form the $7$-simplex with center $D_{1}$. By
corollary \ref{Coro-Skew8-cornered}, the center of $7$-simplex with center
$D_{1}\ $gives the center $D$ of an $A_{3}^{8}(1)$\textbf{-}polytope in
$4_{21}$by $3D+4K_{S_{8}}=D_{1}$. Thus $\left(  D_{1}-4K_{S_{8}}\right)
/3\ $is an integral class in $Pic\ S_{8}$, and we can define a class
$l_{a_{i}}\ $in $Pic\ S_{8}\ $for $a_{i}\ $as%
\[
l_{a_{i}}:=\frac{\left(  D_{1}-K_{S_{8}}\right)  }{3}-a_{i}\ \text{for
}i=1,2,3,4\text{.}%
\]
The classes $l_{a_{i}}$'s are lines because $l_{a_{i}}\cdot K_{S_{8}}=-1\ $and
$l_{a_{i}}^{2}=-1$. Furthermore, we have $l_{a_{i}}\cdot l_{a_{j}}%
=1+a_{i}\cdot a_{j}=1\ $for $i\not =j$ and
\[
l_{a_{1}}+l_{a_{2}}+l_{a_{3}}+l_{a_{4}}=\frac{4\left(  D_{1}-K_{S_{8}}\right)
}{3}-\left(  a_{1}+a_{2}+a_{3}+a_{4}\right)  =D.
\]
Therefore, $l_{a_{i}}$'s\ are four vertices of an uncornered $A_{3}^{8}%
(1)$\textbf{-}polytope in $4_{21}\ $with the center $D$. This gives
theorem.$\ \ \ \ \ \ \ \ \ \ \ \ \ {\LARGE \blacksquare}$

\bigskip

By above proposition and theorem, we have the following configuration of
uncornered $A_{3}^{8}(1)$\textbf{-}polytopes in $4_{21}$.

\bigskip

\begin{theorem}
Let $D\ $be the center of an uncornered $A_{3}^{8}(1)$\textbf{-}polytope in
$4_{21}$ and $D_{1}\ $be the corresponding skew $8$-lines in $4_{21}$. The
configuration of the $A_{3}^{8}(1)$\textbf{-}polytopes in $4_{21}\ $equal the
family of the subset consisting of $4$-skew edges in the $7$-simplex with
center $D_{1}\ $in $4_{21}$.
\end{theorem}

\bigskip

\textbf{(3) Construction of cornered }$A_{3}^{8}(1)$\textbf{-polytopes from
}$A_{2}^{8}(1)$\textbf{-polytopes in }$4_{21}$\textbf{\ }

\bigskip

Recall that each$\ A_{2}^{8}(1)$\textbf{-}polytope in $4_{21}$ is corresponded
to a skew $2$-lines, and the lines consisting the $2$-simplex are in the
vertex figure of each lines in the skew $2$-lines. For example, we have three
lines $l_{1}$, $l_{2}\ $and $l_{3}\ $making an $A_{2}^{8}(1)$\textbf{-}%
polytope, and $l_{a}^{123}+l_{b}^{123}\ $is the corresponding skew
$2$-lines\ with $l_{1}+l_{2}+l_{3}+K_{S_{8}}=l_{a}^{123}+l_{b}^{123}$. Since
$l_{1}$, $l_{2}\ $and $l_{3}\ $are in the vertex figure of $l_{a}^{123}$, the
$A_{2}^{8}(1)$\textbf{-}polytope given by $l_{1}$, $l_{2}\ $and $l_{3}$can be
extended to exactly one cornered $A_{3}^{8}(1)$\textbf{-}polytope in the
vertex figure of $l_{a}^{123}$. Similarly, we can obtain a unique cornered
$A_{3}^{8}(1)$\textbf{-}polytope containing $l_{1}$, $l_{2}\ $and $l_{3}\ $in
the vertex figure of $l_{b}^{123}$. Therefrom, we study the relationship
between above $2$-simplex and the corresponding skew $2$-lines along the
configuration of cornered $A_{3}^{8}(1)$\textbf{-}polytopes.

One of the important feature of a cornered $A_{3}^{8}(1)$\textbf{-}polytope,
there is a line whose vertex figure contains the $3$-simplex. For a line in
$L_{8}\ $and the vertex figure, we consider the blow down map $\pi_{l}%
^{8}:S_{8}\rightarrow S_{7}\ $and observe that the Gieser transform can be
defined on each vertex figure in $4_{21}$. This gives the following definition.

\bigskip

\begin{definition}
Let $l\ $be a line in $L_{8}\ $and $l^{\prime}\ $be a line in the vertex
figure of $l$. The line $G_{l}(l^{\prime})\ $defined as
\[
G_{l}(l^{\prime}):=\pi_{l}^{8\ast}\left(  G\left(  \pi_{l\ast}^{8}(l^{\prime
})\right)  \right)  =-(K_{S_{8}}-l)-l^{\prime}\text{,}%
\]
is called the Gieser transform of $l^{\prime}\ $for $l\ $in $4_{21}$.
\end{definition}

\bigskip

\textbf{Remark:\ }(1) It is easy to see that $G_{l}(l^{\prime})\ $is in the
vertex figure of $l$, and the definition implies $G_{l}(l^{\prime})\cdot
l^{\prime}=2$.

(2) Because $l\cdot l^{\prime}=0$,\ we have two roots $\pm\left(  l-l^{\prime
}\right)  ,\ $and one of the root $l-l^{\prime}\ $determines the line
$G_{l}(l^{\prime})\ $as
\[
G_{l}(l^{\prime})=-K_{S_{8}}+\left(  l-l^{\prime}\right)  \text{.}%
\]
and the other root determine a line $G_{l^{\prime}}(l)$.

\bigskip

The Gieser transform of $l^{\prime}\ $for $l\ $in $4_{21}\ $and the Bertini
transform $B\ $on $4_{21}\ $are related as follows.

\begin{lemma}
\label{Leema-BandG}Let $l\ $and $l^{\prime}\ $be lines in $L_{8}\ $with
$l\cdot l^{\prime}=0$. Then we have (1) $G_{l}(l^{\prime})=B(G_{l^{\prime}%
}(l))\ $and (2)$\ B(l^{\prime})=G_{G_{l}(l^{\prime})}(l)$.
\end{lemma}

\textbf{Proof: }(1)\ Observe $G_{l}(l^{\prime})+G_{l^{\prime}}(l)=-2K_{S_{8}}%
$. Thus $G_{l}(l^{\prime})=-2K_{S_{8}}-G_{l^{\prime}}(l)=B(G_{l^{\prime}}(l))$.

(2) Set $G_{l^{\prime}}(l)=\tilde{l}$. Then $G_{l^{\prime}}(\tilde
{l})=G_{l^{\prime}}(G_{l^{\prime}}(l))=l$. Thus\ by (1), $G_{\tilde{l}%
}(l^{\prime})=B(G_{l^{\prime}}(\tilde{l}))\ $can be written as $G_{G_{l}%
(l^{\prime})}(l)=B(l^{\prime})$.

$\ \ {\LARGE \blacksquare}$

\bigskip

\bigskip The role of the Gieser transform in $4_{21}$ on an $A_{2}^{8}%
(1)$\textbf{-}polytope in $4_{21}$ is explained in the following lemma.

\begin{lemma}
\label{Lemma-3Skewgivesaline}Let $l_{1}$, $l_{2}\ $and $l_{3}\ $be three lines
in $L_{8}\ $giving an $A_{2}^{8}(1)$\textbf{-}polytope in $4_{21}\ $,and
$l_{a}^{123}+l_{b}^{123}\ $is the corresponding skew $2$-lines.

(1)\ The line $G_{l_{a}^{123}}(l_{b}^{123})\ $can be written as $G_{l_{a}%
^{123}}(l_{b}^{123}):=2\left(  l_{a}^{123}-K_{S_{8}}\right)  -D\ $where
$D=(l_{1}+l_{2}+l_{3})$.

(2) $l_{1}$, $l_{2},l_{3}\ $and $G_{l_{a}^{123}}(l_{b}^{123})\ $(resp.
$G_{l_{b}^{123}}(l_{a}^{123})$)\ form a cornered $A_{3}^{8}(1)$\textbf{-}%
polytope in the vertex figure of $l_{a}^{123}\ $(resp. $l_{b}^{123}$).

(3) $G_{l_{a}^{123}}(l_{b}^{123})\ $and $G_{l_{b}^{123}}(l_{a}^{123})\ $are
transformed to each other by the Bertini transform $B$.
\end{lemma}

\textbf{Proof: }(1)\ It is directly given by $(l_{1}+l_{2}+l_{3})+K_{S_{8}%
}=l_{a}^{123}+l_{b}^{123}$.

(2)We just need to check $l_{i}\cdot G_{l_{a}^{123}}(l_{b}^{123})=l_{i}%
\cdot\left(  -K_{S_{8}}+l_{a}^{123}-l_{b}^{123}\right)  =1$ for $i=1,2,3$.

(3) This is given by lemma \ref{Leema-BandG}.

$\ \ {\LARGE \blacksquare}$

\bigskip

\textbf{Remark :} The definition of $G_{l_{a}^{123}}(l_{b}^{123})\ $can be
understood another way. As we consider the center of $A_{2}^{8}(1)$%
\textbf{-}polytope $(l_{1}+l_{2}+l_{3})$, if line $l_{a}^{123}\ $is a line in
$L_{8}$ whose vertex figure contains $l_{1},l_{2}\ $and $l_{3}$, then we can
obtain a line $l^{\prime}\ $given as
\[
l^{\prime}=2\left(  l_{a}^{123}-K_{S_{8}}\right)  -(l_{1}+l_{2}+l_{3}).
\]
And from the above lemma-(1), $G_{l_{a}^{123}}(l_{b}^{123})\ $is $l^{\prime}$.

\bigskip

Now, we consider a cornered $A_{3}^{8}(1)$\textbf{-}polytope in $4_{21}%
\ $where its vertices $l_{1}$, $l_{2}$, $l_{3}\ $and $l_{4}\ $are in the
vertex figure of a line $l$ . From above lemma\ref{Lemma-3Skewgivesaline}, we
have
\[
l_{1}+l_{2}+l_{3}+K_{S_{8}}=l+G_{l}(l_{4}).
\]
Similarly, for each three lines $\{l_{1},l_{2},l_{3},l_{4}\}\backslash
\{l_{i}\},\ i=1,2,3,4$, we get the corresponding $2$-skew lines $l+G_{l}%
(l_{i})$. Here we have another cornered $A_{3}^{8}(1)$\textbf{-}polytope in
$N_{0}(l,S_{8})\ $given by $\{G_{l}(l_{1}),G_{l}(l_{2}),G_{l}(l_{3}%
),G_{l}(l_{4})\}$. This gives the following Theorem.

\bigskip

\begin{theorem}
$\ $Let $l_{1}$, $l_{2}$, $l_{3}\ $and $l_{4}$ be lines in $L_{8}$ which form
an $A_{3}^{8}(1)$\textbf{-}polytope in $4_{21}\ $where$\ $it is in the vertex
figure of a line $l$ . Then $G_{l}\left(  l_{1}\right)  $, $G_{l}\left(
l_{2}\right)  $, $G_{l}(l_{3})\ $and $G_{l}(l_{4})$ form an $A_{3}^{8}%
(1)$\textbf{-}polytope in $4_{21}$.
\end{theorem}

\textbf{Remark:} The $A_{3}^{8}(1)$\textbf{-}polytope\ given by $\{G_{l}%
\left(  l_{1}\right)  ,G_{l}\left(  l_{2}\right)  ,G_{l}(l_{3}),G_{l}%
(l_{4})\}\ $is called the \textit{Gieser dual} of the $A_{3}^{8}(1)$%
\textbf{-}polytope given by $\{l_{1},l_{2},l_{3},l_{4}\}$.

$\bigskip\ $

\subsection{Inscribed $1$-degree $4,5,6\ $and $7$-simplexes in $4_{21}$}

\bigskip

In $4_{21}$, $1$-degree inscribed $m$-simplexes, $A_{m}^{8}(1)$\textbf{-}%
polytopes,exist for $m\leq7$.

The center of each $A_{m}^{8}(1)$\textbf{-}polytope represents a divisor
$l_{1}+...+l_{m+1}$ where $l_{1},...,l_{m+1}$ lines with $l_{i}\cdot l_{j}=1,$
$i\not =j$, which satisfies $\left(  l_{1}+...+l_{m+1}\right)  ^{2}=\left(
m+1\right)  ^{2}-2\left(  m+1\right)  $ and $\left(  l_{1}+...+l_{m+1}\right)
\cdot K_{S_{r}}=-\left(  m+1\right)  $. In particular, for $r=8$, lines
$l\ $in $L_{8}\ $corresponds to roots $l+K_{S_{8}}$, and the divisor class
$l_{1}+...+l_{m+1}$ is transformed to a divisor class $D_{1}:=\left(  \left(
l_{1}+K_{S_{8}}\right)  +...+\left(  l_{m+1}+K_{S_{8}}\right)  \right)  $ with
$D_{1}^{2}=-2\left(  m+1\right)  $ and $D_{1}\cdot K_{S_{8}}=0$. Since
$\left(  l_{i}+K_{S_{8}}\right)  \cdot\left(  l_{j}+K_{S_{8}}\right)
=l_{i}\cdot l_{j}-1$,$\ D_{1}$ is a sum of perpendicular roots. Because the
space of roots in $Pic\ S_{8}\ $is $8$-dimensional, the maximal number of $m$
is $7$. Since cases with $m\leq3$ are discussed the above, $m=4,5,6\ $and $7$
are the remained cases.

According to chapter $4$-\cite{Conway-Sloane}, the set of divisor $D$ with
with $D^{2}=-2\left(  m+1\right)  $ and $D\cdot K_{S_{8}}=0$ corresponds to
the set of vectors with norm $2\left(  m+1\right)  $ in $E_{8}$-lattices.
Here, each set of centers of $A_{m}^{8}(1)$-polytopes $4\leq m\leq7$ in
$4_{21}$ is a subset of the set of vectors with norm $10$, $12,14\ $and $16$
in $E_{8}$-lattice, respectively, and the numbers of the sets of vectors with
norm $10$, $12,14\ $and $16$ in $E_{8}$-lattice are $30240$, $60480,82560\ $%
and $140400$, respectively. Moreover, each set of vertices with norm $10$ and
$12$ which contains divisors $l_{1}+...+l_{m+1}+(m+1)K_{S_{8}}$ for $m=4,5$
are transitively acted by the Weyl group $W(S_{8})\ $(chapter $4$%
-\cite{Conway-Sloane}). Since there is only one orbit of $E_{8}\ $action on
each set of vectors with norm $10\ $and $12$, the numbers of the centers of
$A_{4}^{8}(1)$-polytopes and $A_{5}^{8}(1)$-polytopes are $30240$ and $60480$,
respectively. Here we note the number of skew $3$-lines in $S_{8}$ is also
$60480$.

\bigskip

\subsubsection{Rulings and inscribed $1$-degree $4$-simplexes}

\bigskip

As the inscribed $1$-degree $4$-simplexes, $A_{4}^{8}(1)$\textbf{-}polytopes
exist only for $r=8$, there is no line whose vertex figure contains an
$A_{4}^{8}(1)$\textbf{-}polytope. Thus every$\ A_{4}^{8}(1)$\textbf{-}%
polytope\ is uncornered. By adding a proper line to a cornered $A_{3}^{8}%
(1)$\textbf{-}polytope, we can obtain an $A_{4}^{8}(1)$\textbf{-}polytope. But
the choice of new line is too arbitrary to describe the configuration of
$A_{4}^{8}(1)$\textbf{-}polytopes. Thus we begin with $A_{4}^{8}(1)$%
\textbf{-}polytopes constructed from uncornered $A_{3}^{8}(1)$\textbf{-}%
polytopes in $4_{21}$. This leads us the configuration of $A_{4}^{8}%
(1)$\textbf{-}polytopes in $4_{21}$.

Let $l_{1}$, $l_{2}$, $l_{3}\ $and $l_{4}$ be lines in $L_{8}$ which form an
uncornered $A_{3}^{8}(1)$\textbf{-}polytope\ in $4_{21}$. We consider the skew
$2$-lines $l_{a}^{123}+l_{b}^{123}\ $corresponding to $l_{1}+l_{2}+l_{3}$,
namely, $l_{1}+l_{2}+l_{3}+K_{S_{8}}=l_{a}^{123}+l_{b}^{123}$. By lemma
\ref{Lemma-3Skewgivesaline}, we have a cornered $A_{3}^{8}(1)$\textbf{-}%
polytope in $4_{21}\ $given by $\{l_{1},l_{2},l_{3},G_{l_{a}^{123}}%
(l_{b}^{123})\}$. In fact we obtain an $A_{4}^{8}(1)$\textbf{-}polytope\ in
$4_{21}\ $as follows.

\bigskip

\begin{theorem}
\label{Thm-uncorner-4to5}Suppose $l_{1}$, $l_{2}$, $l_{3}\ $and $l_{4}$, lines
in $L_{8}$ give an uncornered $A_{3}^{8}(1)$\textbf{-}polytope in $4_{21}$,
and $l_{a}^{123}+l_{b}^{123}\ $is the skew $2\ $-lines corresponding to
$l_{1}+l_{2}+l_{3}$. The five lines $l_{1}$, $l_{2}$, $l_{3},G_{l_{a}^{123}%
}(l_{b}^{123})\ $and $l_{4}\ $form an $A_{4}^{8}(1)$\textbf{-}polytope in
$4_{21}$.
\end{theorem}

\textbf{Proof: }We only need to check $G_{l_{a}^{123}}(l_{b}^{123})\cdot
l_{4}=1$, and we need to get $l_{a}^{123}\cdot l_{4}\ $and $l_{b}^{123}\cdot
l_{4}\ $at the first hand.

We consider another skew $2\ $-lines $l_{a}^{124}+l_{b}^{124}\ $given by
$l_{1}+l_{2}+l_{4}+K_{S_{8}}=l_{a}^{124}+l_{b}^{124}$. Since the
$3$-simplex\ is uncornered, by theorem \ref{Thm-uncorner-skewedge}
$l_{a}+l_{b}\ $and $l_{\alpha}+l_{\beta}\ $are skew edges. Thus we have
$l_{a}^{123}\cdot\left(  l_{b}^{124}+l_{b}^{124}\right)  =l_{b}^{123}%
\cdot\left(  l_{a}^{124}+l_{b}^{124}\right)  =0$. Furthermore, $l_{a}%
^{123}\cdot l_{4}=1\ $because
\[
0=l_{a}^{123}\cdot\left(  l_{a}^{124}+l_{b}^{124}\right)  =l_{a}^{123}%
\cdot\left(  l_{1}+l_{2}+l_{4}+K_{S_{8}}\right)  =l_{a}^{123}\cdot
l_{4}-1\text{,}%
\]
and the same reason $l_{b}^{123}\cdot l_{4}=1$.

Now, we have $\ $%
\[
G_{l_{a}^{123}}(l_{b}^{123})\cdot l_{4}=\left(  -K_{S_{8}}+l_{a}^{123}%
-l_{b}^{123}\right)  \cdot l_{4}=1,
\]
and $\{l_{1},l_{2},l_{3},G_{l_{a}^{123}}(l_{b}^{123}),l_{4}\}\ $gives an
inscribed $1$-degree $4$-simplex in $4_{21}$.$\ \ \ $

$\ {\LARGE \blacksquare}$

\bigskip

\textbf{Remark : }This theorem implies that any line $l\ $where $\{l_{1}%
,l_{2},l_{3},l\}\ $gives an inscribed uncornered $1$-degree $3$-simplex in
$4_{21}$, $G_{l_{a}^{123}}(l_{b}^{123})\cdot l=G_{l_{b}^{123}}(l_{a}%
^{123})\cdot l=1$.

\bigskip

In the theorem \ref{Thm-uncorner-4to5}, the $A_{4}^{8}(1)$\textbf{-}polytope
in $4_{21}\ $given by $\{l_{1},l_{2},l_{3},G_{l_{a}^{123}}(l_{b}^{123}%
),l_{4}\}\ $contains a cornered $A_{3}^{8}(1)$\textbf{-}polytope consisting of
$\{l_{1},l_{2},l_{3},G_{l_{a}^{123}}(l_{b}^{123})\}$. In fact this is the only
possible $1$-degree $3$-simplex\ in the $4$-simplex\ because of the
proposition \ref{Prop-4cornerTo3+1Uncornered}. It turns out this is true for
all the $A_{4}^{8}(1)$\textbf{-}polytopes in $4_{21}\ $as follows.

\bigskip

\begin{theorem}
\label{Thm-4sym-with cornered}Let $\{l_{i}\ 1\leq i\leq5\}\ $be the set of
lines in $L_{8}\ $giving an $A_{4}^{8}(1)$\textbf{-}polytope in $4_{21}$.
There is a unique subset of $\{l_{i}\ 1\leq i\leq5\}$ which gives a cornered
$A_{3}^{8}(1)$\textbf{-}polytope.
\end{theorem}

\textbf{Proof:} \textbf{(Existence) }Suppose all the $A_{3}^{8}(1)$%
\textbf{-}polytopes in $\{l_{i}\ 1\leq i\leq5\}\ $are uncornered. Then we
consider an $A_{2}^{8}(1)$\textbf{-}polytope\ given by $\{l_{1},l_{2}%
,l_{3}\}\ $and the corresponding skew $2$-lines $l_{a}^{123}+l_{b}^{123}$. By
lemma \ref{Lemma-3Skewgivesaline}\ and the remark of the theorem
\ref{Thm-uncorner-4to5}, we have $G_{l_{a}^{123}}(l_{b}^{123})\cdot
l_{4}=G_{l_{a}^{123}}(l_{b}^{123})\cdot l_{5}=1\ $because $\{l_{1},l_{2}%
,l_{3},l_{4}\}\ $and $\{l_{1},l_{2},l_{3},l_{5}\}\ $are uncornered $A_{3}%
^{8}(1)$\textbf{-}polytopes. Thus $\{l_{1},l_{2},l_{3},l_{4},l_{5}%
,G_{l_{a}^{123}}(l_{b}^{123})\}\ $gives an $A_{5}^{8}(1)$\textbf{-}polytope in
$4_{21}$. Similarly, for a skew $2$-lines $l_{a}^{124}+l_{b}^{124}%
\ $corresponding to the inscribed $1$-degree $2$-simplex of $\{l_{1}%
,l_{2},l_{4}\}$, we obtain another$\ A_{5}^{8}(1)$\textbf{-}polytope
$\{l_{1},l_{2},l_{3},l_{4},l_{5},G_{l_{a}^{124}}(l_{a}^{124})\}$ in $4_{21}$.
Furthermore, we observe that $G_{l_{a}^{123}}(l_{b}^{123})\cdot G_{l_{a}%
^{124}}(l_{b}^{124})=1\ $verified as%
\[
G_{l_{a}^{123}}(l_{b}^{123})\cdot G_{l_{a}^{124}}(l_{b}^{124})=\left(
-K_{S_{8}}+l_{a}^{123}-l_{b}^{123}\right)  \cdot\left(  -K_{S_{8}}+l_{a}%
^{124}-l_{b}^{124}\right)  =1.
\]
Here $l_{a}^{123}\cdot l_{a}^{124}=l_{a}^{123}\cdot l_{b}^{124}=l_{b}%
^{123}\cdot l_{a}^{124}=l_{b}^{123}\cdot l_{b}^{124}=0\ $because $l_{a}%
^{123}+l_{b}^{123}\ $and $l_{a}^{124}+l_{b}^{124}\ $are two skew edges given
from $\{l_{1},l_{2},l_{3},l_{4}\}$ by the proposition \ref{Prop-skewedges}%
.\ Thus any three lines from $\{l_{1},l_{2},l_{3},l_{4}\}\ $gives a line so
that the
\[
\{l_{1},l_{2},l_{3},l_{4},l_{5},G_{l_{a}^{123}}(l_{b}^{123}),G_{l_{a}^{124}%
}(l_{b}^{124}),G_{l_{a}^{134}}(l_{b}^{134}),G_{l_{a}^{234}}(l_{b}^{234})\}
\]
$\ $makes an $A_{8}^{8}(1)$\textbf{-}polytope\ in $4_{21}$. But this
$8$-simplex in $4_{21}$ may not exist in $4_{21}$. Thus we conclude that there
exists a subset in $\{l_{i}\ 1\leq i\leq5\}$ which gives an $A_{3}^{8}%
(1)$\textbf{-}polytope.

\textbf{(Uniqueness)}\ Assume $\{l_{1},l_{2},l_{3},l_{4}\}\ $is a cornered
$A_{3}^{8}(1)$\textbf{-}polytope from $\{l_{i}\ 1\leq i\leq5\}$. There are
four other possible choices of $A_{3}^{8}(1)$\textbf{-}polytopes in
$\{l_{i}\ 1\leq i\leq5\}$ which are obtained by replacing one of the lines in
$\{l_{1},l_{2},l_{3},l_{4}\}\ $by $l_{5}$. But by proposition
\ref{Prop-4cornerTo3+1Uncornered}, none of these four $3$-simplexes can be
cornered. Thus there is only one cornered $A_{3}^{8}(1)$\textbf{-}%
polytope\ from $\{l_{i}\ 1\leq i\leq5\}$.

$\ \ {\LARGE \blacksquare}$

\textbf{ }

Recall the number of center of $A_{4}^{8}(1)$\textbf{-}polytopes in $4_{21}%
\ $is $30240$. We can observe this number equals to the number of ordered
pairs of lines in $L_{8}\ $with $1$-intersection, namely, $\left\vert
L_{8}\right\vert \left\vert N_{1}(l,S_{8})\right\vert =240\times
126=240\times(5^{3}+1^{3})$. In the following theorem, we show that the
coincidence of these two number leads us the configuration of $A_{4}^{8}%
(1)$\textbf{-}polytopes in $4_{21}$.

\bigskip

\begin{theorem}
\label{Thm-Conf-4in8}Let $D\ $be a center of an $A_{4}^{8}(1)$\textbf{-}%
polytope in $4_{21}$. There is unique line $l_{D}^{S_{8}}\ $in $L_{8}\ $and
the center $A_{D}\ $of an cornered $A_{3}^{8}(1)$\textbf{-}polytope such that
$D=$ $A_{D}+l_{D}^{S_{8}}$. Furthermore, the set of the centers of $A_{4}%
^{8}(1)$\textbf{-}polytope in $4_{21}\ $is bijective to the following set
$\tilde{F}_{8}\ $defined as
\[
\tilde{F}_{8}:=\{(l_{1},l_{2})\mid l_{1}\text{,}l_{2}\in L_{8}\ \text{with
}l_{1}\cdot l_{2}=1\}\text{.}%
\]

\end{theorem}

\textbf{Proof:} From the above theorem \ref{Thm-4sym-with cornered}, each
center $D\ $of $A_{4}^{8}(1)$\textbf{-}polytopes in $4_{21}\ $can be written
as the sum of a center $A_{1}$ of a cornered $A_{3}^{8}(1)$\textbf{-}polytope
and a line $l_{1}\ $in $L_{8}$. Suppose $A_{2}\ $and $l_{2}\ $are another
center of a cornered $A_{3}^{8}(1)$\textbf{-}polytope\ and a line in $L_{8}$
with $D=A_{1}+l_{1}=A_{2}+l_{2}$, $(A_{1},l_{1})\not =(A_{2},l_{2})$.

Since $A_{1}\ $is the center of a cornered $A_{3}^{8}(1)$\textbf{-}polytope,
there is a line $l\ $where its vertex figure contains the $3$-simplex. Thus
$A_{1}\cdot l=0$,$\ $and $D\cdot l=1\ $because $l_{1}\cdot l=l_{1}\cdot\left(
K_{S_{8}}+\frac{A_{1}}{2}\right)  =1$. This implies that $l\cdot A_{2}+l\cdot
l_{2}=1$, and we consider the following cases for $l\cdot l_{2}$. Since
$l_{1}\cdot A_{2}\geq-1$, we only consider $-1\leq l\cdot l_{2}\leq2$.

(a) If $l\cdot l_{2}=-1$, then $l\cdot A_{2}=2$.\ Since $l=l_{2}$, we have
$l_{2}\cdot A_{2}=l\cdot A_{2}=2$. But by the definition of the $A_{4}^{8}%
(1)$\textbf{-}polytope, $l_{2}\cdot A_{2}\ $must be $4$. Thus $l\cdot
l_{2}\not =-1$.\ (b) If $l\cdot l_{2}=0$, then $l\cdot A_{2}=1$. Here $l\ $can
not be in $B$ the set of the lines giving the cornered $A_{3}^{8}%
(1)$\textbf{-}polytope\ with center $A_{2}$, because in this case $l\cdot$
$A_{2}\ $must be $2\ $instead of $1\ $by the definition of the $A_{4}^{8}%
(1)$\textbf{-}polytope. Now since $l\ $is not in $B$, $l\ $intersects
positively with only one line $l_{B}$ in $B$,$\ $and the other lines in
$B\ $are in the vertex figure of $l$. Furthermore, $B\backslash\{l_{B}\}\ $and
$l_{2}\ $form an $A_{3}^{8}(1)$\textbf{-}polytope which is cornered by the
line $l$. But by theorem \ref{Thm-4sym-with cornered}, this $A_{3}^{8}%
(1)$\textbf{-}polytope can not be cornered. Thus $l\cdot l_{2}\not =0$.\ (c)
If $l\cdot l_{2}=1$, then $l\cdot A_{2}=0$. Again, $l\ $is not one of the line
in $B$, and thus $l\ $has $0$-intersection to each line in $B$. This implies
that $l\ $and $A_{2}\ $are corresponded by corollary \ref{CoroCornered-line}.
Therefore, we have $A_{1}=A_{2}\ $with is a contradiction. Thus $l\cdot
l_{2}\not =1$.\ (d) If $l\cdot l_{2}=2$, then $l\cdot A_{2}=-1$. Here
$l\ $must be one of the line in $B$. But in that case $l\cdot A_{2}%
=2\ $instead of $-1$. Thus $l\cdot l_{2}\not =-1$.

From above cases, we conclude the choice of center $A_{1}\ $and a line
$l_{1}\ $for $D=A_{1}+l_{1}\ $is unique.

By corollary\ref{CoroCornered-line}, we get a line $l^{\prime}\ $from the
center $A_{1}$, and $l^{\prime}\ $has $l^{\prime}\cdot l_{1}=1$. Thus the
ordered pair $(l^{\prime},l_{1})\ $in $\tilde{F}_{8}$. Obviously, each pair of
lines in $\tilde{F}_{8}\ $gives a center of an $A_{4}^{8}(1)$\textbf{-}%
polytope in $4_{21}$.

Thus we have theorem.

$\ \ {\LARGE \blacksquare}$

\bigskip

By the above theorem, we have the following theorem for configuration of
$A_{4}^{8}(1)$\textbf{-}polytopes in $4_{21}$

\begin{theorem}
\label{Thm-Conf-4} Suppose $D\ $is a center of an $A_{4}^{8}(1)$%
\textbf{-}polytope in $4_{21}$. All the $A_{4}^{8}(1)$\textbf{-}polytopes with
center $D$ in $4_{21}\ $share a common line $l$, and $D-l\ $is the common
center for the unique cornered $A_{3}^{8}(1)$\textbf{-}polytope in
each$\ A_{4}^{8}(1)$\textbf{-}polytope.
\end{theorem}

\bigskip

\subsubsection{Skew $3$-lines and inscribed $1$-degree\ $5$-simplexes}

\bigskip

Recall that the number of center of inscribed $1$-degree $5$-simplexes,
$A_{5}^{8}(1)$\textbf{-}polytopes in $4_{21}\ $equals to the number of skew
$3$-lines in $4_{21}$ where each skew $3$-lines can be written as the sum of
unique three lines in $L_{8}$. While the $A_{5}^{8}(1)$\textbf{-}polytopes do
not imply the uniqueness issues directly, the uniqueness appearing in the skew
$3$-lines leads us the following lemma, and therefrom, we get the
configuration of $A_{5}^{8}(1)$\textbf{-}polytopes in $4_{21}$.

\bigskip

\begin{lemma}
\label{Lemma-3-lemma}\ For each six lines in $L_{8}$ consisting of an
inscribed $1$-degree $5$-simplex in $4_{21}$, there exist three cornered
$A_{3}^{8}(1)$\textbf{-}polytope in it, and the choice is unique. Furthermore,
the given six lines can be labelled as $l_{i}\ 1\leq i\leq6$\ \ where
$\{l_{1},l_{2},l_{3},l_{4}\}$,$\{l_{1},l_{2},l_{5},l_{6}\}\ $and
$\{l_{3},l_{4},l_{5},l_{6}\}\ $are the cornered $A_{5}^{8}(1)$\textbf{-}polytopes.
\end{lemma}

\textbf{Proof :} Recall that for each five lines from the given six lines,
there exists an cornered $A_{3}^{8}(1)$\textbf{-}polytope by the theorem
\ref{Thm-4sym-with cornered}. We choose four lines from the six lines which
give a cornered $A_{3}^{8}(1)$\textbf{-}polytope, and label them as
$l_{1},l_{2},l_{3}\ $and$\ l_{4}$. By the proposition
\ref{Prop-4cornerTo3+1Uncornered}, the lines from another cornered $A_{3}%
^{8}(1)$\textbf{-}polytope and $\{l_{1},l_{2},l_{3},l_{4}\}\ $share $1\ $or
$2\ $lines. Since there are only six lines, they share $2$-lines, and we can
relabel six lines so that $\{l_{1},l_{2},l_{3},l_{4}\}\ $and $\{l_{1}%
,l_{2},l_{5},l_{6}\}\ $produce cornered $A_{3}^{8}(1)$\textbf{-}polytopes. By
applying theorem \ref{Thm-4sym-with cornered}\ and the proposition
\ref{Prop-4cornerTo3+1Uncornered} to each $\{l_{1},l_{2},l_{3},l_{4}%
,l_{5},l_{6}\}\backslash\{l_{i}\}\ 1\leq i\leq6$, we deduce that there is one
more cornered $A_{3}^{8}(1)$\textbf{-}polytope\ given by $\{l_{3},l_{4}%
,l_{5},l_{6}\}\ $which also share two lines with $\{l_{1},l_{2},l_{3}%
,l_{4}\}\ $and $\{l_{1},l_{2},l_{5},l_{6}\}$. This gives the lemma.

$\ {\LARGE \blacksquare}$

\bigskip

\textbf{Remark:} According to the lemma\ref{Lemma-3-lemma}, each set
$A\ $of\ six lines in $L_{8}$ consisting of an $A_{5}^{8}(1)$\textbf{-}%
polytope in $4_{21}$, there are unique three disjoint subsets $A_{i}%
\ i=1,2,3\ $of $A\ $consisting of two lines such that each $A-A_{i}$
$i=1,2,3\ $produces a cornered $A_{3}^{8}(1)$\textbf{-}polytope.

\bigskip

\begin{theorem}
\label{Thm-conf5}Let $D\ $be the center of an $A_{5}^{8}(1)$\textbf{-}polytope
in $4_{21}$. The center $D\ $is corresponded to a skew $3$-lines $l_{a}%
+l_{b}+l_{c}\ $by $D+3K_{S_{8}}=l_{a}+l_{b}+l_{c}$. Furthermore, each
$A_{5}^{8}(1)$\textbf{-}polytope in $4_{21}\ $with center $D\ $has
three\ uniquely determined cornered $A_{3}^{8}(1)$\textbf{-}polytopes where
they are in the vertex figure of $l_{a}$,\ $l_{b}\ $and $l_{c}\ $respectively.
\end{theorem}

\textbf{Proof :}\ Since the center $D\ $satisfying $D^{2}=24$, $D\cdot
K_{S_{8}}=-6$, the new divisor class $D+3K_{S_{8}}\ $also satisfies $\left(
D+3K_{S_{8}}\right)  ^{2}=-3$,$\left(  D+3K_{S_{8}}\right)  \cdot K_{S_{8}%
}=-3$, and it is a skew $3$-line. Since a skew $3$-line is the sum of unique
three disjoint lines, we have $D+3K_{S_{8}}=l_{a}+l_{b}+l_{c}\ $where $l_{a}%
$,$\ l_{b}\ $and $l_{c}\ $are disjoint lines in $L_{8}$.

On the other hand, by above lemma \ref{Lemma-3-lemma}, we label the given six
lines as $l_{i}\ 1\leq i\leq6$\ so that $\{l_{1},l_{2},l_{3},l_{4}\}$%
,$\{l_{1},l_{2},l_{5},l_{6}\}\ $and $\{l_{3},l_{4},l_{5},l_{6}\}\ $are
cornered $A_{3}^{8}(1)$\textbf{-}polytopes. By corollary
\ref{CoroCornered-line}, these cornered $A_{3}^{8}(1)$\textbf{-}%
polytopes\ produce lines $l^{a},l^{b}\ $and $l^{c}$ in $L_{8}$. Here, the
intersection $l^{a}\cdot l^{b}\ $is%
\[
\left(  \frac{l_{1}+l_{2}+l_{3}+l_{4}}{2}+K_{S_{8}}\right)  \cdot\left(
\frac{l_{1}+l_{2}+l_{5}+l_{6}}{2}+K_{S_{8}}\right)  =0\text{,}%
\]
and similarly $l^{a}$,$\ l^{b}\ $and $l^{c}\ $are disjoint each other.
Furthermore, we have%
\begin{align*}
l^{a}+l^{b}+l^{c}  &  =(l_{1}+l_{2}+l_{3}+l_{4}+l_{5}+l_{6})+3K_{S_{8}}\\
&  =D+3K_{S_{8}}=l_{a}+l_{b}+l_{c}\text{.}%
\end{align*}
Thus $l^{a}$,$\ l^{b}\ $and $l^{c}\ $form the skew $3$-line $l_{a}+l_{b}%
+l_{c}$,$\ $and therefore $\{l^{a},l^{b},l^{c}\}=\{l_{a},l_{b},l_{c}\}$. This
gives the theorem.

$\ {\LARGE \blacksquare}$

\bigskip

\subsubsection{Fano planes and inscribed $1$-degree\ $6$-simplexes}

\bigskip

The set of center of inscribed $1$-degree $6$-simplexes, $A_{6}^{8}%
(1)$\textbf{-}polytopes in $4_{21}\ $is a subset of the set of vectors with
norm $14$ in $E_{8}$-lattice which has more than one $E_{8}$-orbit. But we
figure out that the configuration of $A_{6}^{8}(1)$\textbf{-}polytopes in
$4_{21}\ $is somewhat similar to that of $A_{5}^{8}(1)$\textbf{-}polytopes.
Here each center of $A_{6}^{8}(1)$\textbf{-}polytope gives a skew $7$-lines in
$4_{21}\ $where the choice of lines in it is unique, and the uniqueness is the
key to study the configuration of $A_{6}^{8}(1)$\textbf{-}polytopes in
$4_{21}$.

\bigskip

\begin{lemma}
\label{Lemma-7lemma}For each seven lines in $L_{8}$ consisting of an
$A_{6}^{8}(1)$\textbf{-}polytope in $4_{21}$, there exist seven cornered
$A_{3}^{8}(1)$\textbf{-}polytopes in it, and the choice is unique.
Furthermore, the given seven lines can be labelled as $l_{i}\ 1\leq i\leq
7$\ \ where
\begin{align*}
&  \{l_{1},l_{2},l_{3},l_{4}\},\{l_{1},l_{2},l_{5},l_{6}\},\{l_{3},l_{4}%
,l_{5},l_{6}\},\\
&  \{l_{1},l_{3},l_{5},l_{7}\},\{l_{2},l_{4},l_{5},l_{7}\},\{l_{1},l_{4}%
,l_{6},l_{7}\}\ \\
&  \text{and }\{l_{2},l_{3},l_{6},l_{7}\}
\end{align*}
$\ $are the cornered $A_{3}^{8}(1)$\textbf{-}polytopes.
\end{lemma}

\textbf{Proof : }We choose six lines from the given seven lines. By applying
the lemma \ref{Lemma-3-lemma} on the six lines, we can label them as
$l_{i}\ 1\leq i\leq6$\ \ where\ $\{l_{1},l_{2},l_{3},l_{4}\},\{l_{1}%
,l_{2},l_{5},l_{6}\}\ $and $\{l_{3},l_{4},l_{5},l_{6}\}\ $are the three
cornered $A_{3}^{8}(1)$\textbf{-}polytopes in it. Thus the remained one line
is denoted by $l_{7}$.

\ Now, we consider another set of six lines $\{l_{1},l_{2},l_{3},l_{4}%
,l_{5},l_{7}\}$. By lemma \ref{Lemma-3-lemma} and proposition
\ref{Prop-4cornerTo3+1Uncornered}, $\{l_{1},l_{3},l_{5},l_{7}\}$,
$\{l_{1},l_{4},l_{5},l_{7}\},\{l_{2},l_{4},l_{5},l_{7}\}\ $and $\{l_{2}%
,l_{3},l_{5},l_{7}\}\ $are possible subsets of cornered $A_{3}^{8}%
(1)$\textbf{-}polytopes in addition to one given by $\{l_{1},l_{2},l_{3}%
,l_{4}\}$. Here, we observe that even if we exchange $l_{1}\ $(resp. $l_{3}$)
with $l_{2}\ $(resp. $l_{4}$), the $3$-simplexes from $\{l_{1},l_{2}%
,l_{3},l_{4},l_{5},l_{6}\}\ $do not change. Thus, after relabelling if it is
necessary, we can set $\{l_{1},l_{3},l_{5},l_{7}\}\ $and $\{l_{2},l_{4}%
,l_{5},l_{7}\}\ $give new cornered $A_{3}^{8}(1)$\textbf{-}polytopes. Next, we
consider another set of six lines $\{l_{1},l_{2},l_{3},l_{4},l_{6},l_{7}\}$.
Here, by lemma \ref{Lemma-3-lemma} and proposition
\ref{Prop-4cornerTo3+1Uncornered}, $\{l_{1},l_{4},l_{6},l_{7}\}\ $and
$\{l_{2},l_{3},l_{6},l_{7}\}\ $are the only possible subsets of cornered
$A_{3}^{8}(1)$\textbf{-}polytopes from $\{l_{1},l_{2},l_{3},l_{4},l_{6}%
,l_{7}\}$. And one can check there is no more cornered $A_{3}^{8}%
(1)$\textbf{-}polytope from the other choices of six lines from $\{l_{i}%
\ 1\leq i\leq7\ \}$.

By the above algorithm, we show that each seven lines in $L_{8}$ consisting of
an $A_{6}^{8}(1)$\textbf{-}polytope in $4_{21}\ $has seven cornered $A_{3}%
^{8}(1)$\textbf{-}polytopes and the choice is unique. This gives the lemma.

$\ {\LARGE \blacksquare}$

\bigskip

By applying the argument in theorem \ref{Thm-conf5} \ and the lemma
\ref{Lemma-7lemma} to the center of an $A_{6}^{8}(1)$\textbf{-}polytope in
$4_{21}$, we obtain the following theorem. We leave the detail to readers.

\begin{theorem}
\label{Thm-conf6}Let $D\ $be the center of an $A_{6}^{8}(1)$\textbf{-}polytope
in $4_{21}$. The center $D\ $is corresponded to a skew $7$-lines $\sum
_{i=1}^{7}l_{a_{i}}\ $by $2D+7K_{S_{8}}=\sum_{i=1}^{7}l_{a_{i}}$. Furthermore,
each $A_{6}^{8}(1)$\textbf{-}polytope in $4_{21}\ $with center $D\ $has
seven\ uniquely determined cornered $A_{3}^{8}(1)$\textbf{-}polytopes where
they are in the vertex figure of $l_{a_{i}}\ 1\leq i\leq7\ $respectively.
\end{theorem}

\bigskip

In the lemma \ref{Lemma-7lemma}, we observe that when we have four lines for
a$\ $cornered $A_{3}^{8}(1)$\textbf{-}polytope, the remain three lines in the
$A_{6}^{8}(1)$\textbf{-}polytope give an $A_{2}^{8}(1)$\textbf{-}polytope
which is not contained in any$\ $cornered $A_{3}^{8}(1)$\textbf{-}polytope in
the $6$-simplex. The converse is also true by the following proposition.

\begin{proposition}
\label{Prop-7-Remain4-Uncornered}For seven lines in $L_{8}\ $consisting an
$A_{6}^{8}(1)$\textbf{-}polytope in $4_{21}$, if three lines in the
$6$-simplex form an $A_{2}^{8}(1)$\textbf{-}polytope which is not contained in
any$\ $cornered $A_{3}^{8}(1)$\textbf{-}polytope in the $6$-simplex, the
remain four lines in the $7$-simplex give a$\ $cornered $A_{3}^{8}%
(1)$\textbf{-}polytope.
\end{proposition}

\textbf{Proof : }Let $l_{i}\ 1\leq i\leq7\ $be the given seven lines in
$L_{8}\ $consisting an $A_{6}^{8}(1)$\textbf{-}polytope in $4_{21}$. Suppose
$\{l_{5},l_{6},l_{7}\}\ $gives an $A_{2}^{8}(1)$\textbf{-}polytope\ which is
not contained in any$\ $cornered $A_{3}^{8}(1)$\textbf{-}polytope, we want to
show that $\{l_{1},l_{2},l_{3},l_{4}\}\ $gives an cornered$\ A_{3}^{8}%
(1)$\textbf{-}polytope by following the algorithm in the lemma
\ref{Lemma-7lemma}.

First, we claim that the intersection between the subset of lines producing a
cornered $A_{3}^{8}(1)$\textbf{-}polytope and $\{l_{5},l_{6},l_{7}\}\ $must
has two or no lines. By the assumption, the intersection can not be three
lines. Furthermore, if it is a line, the union of the lines for a$\ $cornered
$A_{3}^{8}(1)$\textbf{-}polytope and $\{l_{5},l_{6},l_{7}\}\ $must be set of
six lines and have a$\ $cornered $A_{3}^{8}(1)$\textbf{-}polytope containing
$\{l_{5},l_{6},l_{7}\}$ by lemma \ref{Lemma-7lemma}.

We consider a subset $\{l_{i}\ 2\leq i\leq7\}\ $of the given seven lines. From
the above claim, here are three$\ $cornered $A_{3}^{8}(1)$\textbf{-}polytopes
and each of them contains two lines from $\{l_{5},l_{6},l_{7}\}$.\ Thus
a$\ $cornered $A_{3}^{8}(1)$\textbf{-}polytope\ in $\{l_{i}\ 2\leq i\leq
7\}\ $consists of two lines from $\{l_{2},l_{3},l_{4}\}\ $and two lines from
$\{l_{5},l_{6},l_{7}\}$.Without losing generality, we assume the cornered
$A_{3}^{8}(1)$\textbf{-}polytope is given by $\{l_{2},l_{3},l_{5},l_{6}\}$.
Now, we consider $\{l_{i}\ 1\leq i\leq6\}$. Since $\{l_{2},l_{3},l_{5}%
,l_{6}\}\ $produces a cornered $A_{3}^{8}(1)$\textbf{-}polytope, by the above
claim, we obtain two other cornered $A_{3}^{8}(1)$\textbf{-}polytopes are
given by $\{l_{1},l_{2},l_{5},l_{6}\}\ $and $\{l_{1},l_{2},l_{3},l_{4}\}$.
Thus $\{l_{1},l_{2},l_{3},l_{4}\}\ $must produce a$\ $cornered $A_{3}^{8}%
(1)$\textbf{-}polytope.\ This gives the proposition.

$\ {\LARGE \blacksquare}$

\bigskip

\textbf{Remark : \ }We call the inscribed $1$-degree $2$-simplex in an
$A_{6}^{8}(1)$\textbf{-}polytope\ in the above proposition a \textit{Fano
block }in the\textit{ }$A_{6}^{8}(1)$\textbf{-}polytope which appears in the
following theorem.

\bigskip

From the theorem and the proposition, we have the following theorem.

\begin{theorem}
\label{Thm-Fano}Let $\{l_{i}\ 1\leq i\leq7\}\ $be a set of lines in $L_{8}%
\ $which forms an $A_{6}^{8}(1)$\textbf{-}polytope in $4_{21}$.\ The set of
Fano blocks in the\textit{ }$A_{6}^{8}(1)$\textbf{-}polytope consists of seven
$A_{2}^{8}(1)$\textbf{-}polytopes in the $6$-simplex. Furthermore, the seven
lines and Fano blocks in the $6$-simplex produce a Steiner system
$S(2,3,7)\ $which is known as Fano plane.
\end{theorem}

\bigskip\textbf{Proof :} By proposition \ref{Prop-7-Remain4-Uncornered} and
theorem\ \ref{Thm-conf6}, there are only seven Fano blocks in the given set of
lines. By lemma \ref{Lemma-7lemma}, we can label the seven Fano blocks as%
\begin{align*}
&  \{l_{5},l_{6},l_{7}\},\{l_{3},l_{4},l_{7}\},\{l_{1},l_{2},l_{7}%
\},\{l_{2},l_{4},l_{6}\},\\
&  \{l_{1},l_{3},l_{6}\},\{l_{2},l_{3},l_{5}\}\ \text{and }\{l_{1},l_{4}%
,l_{5}\}\text{.}%
\end{align*}
Furthermore, one can directly check the seven lines $\{l_{i}\ 1\leq
i\leq7\}\ \ $and Fano blocks form a Steiner system $S(2,3,7)$.

\bigskip$\ {\LARGE \blacksquare}$

\subsubsection{Rulings and inscribed $1$-degree\ $7$-simplexes}

Again, the set of center of inscribed $1$-degree $7$-simplexes, $A_{7}^{8}%
(1)$\textbf{-}polytopes, in $4_{21}\ $is a subset of the set of vectors with
norm $16$ in $E_{8}$-lattice which has more than one $E_{8}$-orbit. Here we
find out that the configuration of $A_{7}^{8}(1)$\textbf{-}polytopes in
$4_{21}\ $is close to that of $A_{4}^{8}(1)$\textbf{-}polytopes. Here each
center of $A_{7}^{8}(1)$\textbf{-}polytopes gives a ruling in $4_{21}$.
Furthermore, each center of $A_{7}^{8}(1)$\textbf{-}polytopes can be written
as the sum of two centers of$\ $cornered $A_{3}^{8}(1)$\textbf{-}polytopes
where they are in the vertex figures of the bipolar pair of lines in the
corresponding ruling. The detail is as follows.

\bigskip

\begin{theorem}
Let $D\ $be the center of an $A_{7}^{8}(1)$\textbf{-}polytope in $4_{21}$. The
center $D\ $is corresponded to a ruling $f$ by $D/2+2K_{S_{8}}=f$.
Furthermore, the center $D\ $can be written as a sum of two centers of
cornered $A_{3}^{8}(1)$\textbf{-}polytopes in the $7$-simplex. There are seven
pairs of centers of cornered $A_{3}^{8}(1)$\textbf{-}polytopes whose sum is
$D$, and each pair corresponds the bipolar pairs of lines in the
$7$-crosspolytope given by the ruling $f$.
\end{theorem}

\textbf{Proof :\ }Let $\{l_{i}\ 1\leq i\leq8\}\ $be a set of lines in
$L_{8}\ $producing an $A_{7}^{8}(1)$\textbf{-}polytope in $4_{21}\ $whose
center is $D$.\ We claim that if a subset of four lines in $\{l_{i}\ 1\leq
i\leq8\}\ $gives a cornered $A_{3}^{8}(1)$\textbf{-}polytope with center
$D_{1}$, the complement subset also form another cornered $A_{3}^{8}%
(1)$\textbf{-}polytope\ with center $D_{2}:=D-D_{1}$.

Proof of claim: Without losing generality, suppose $D_{1}\ $is given as
$l_{1}+l_{2}+l_{3}+l_{4}$. If we consider $\{l_{i}\ 1\leq i\leq7\}$,
$\{l_{5},l_{6},l_{7}\}\ $must form a Fano block in $\{l_{i}\ 1\leq i\leq7\}$
by proposition \ref{Prop-7-Remain4-Uncornered} and theorem \ref{Thm-Fano}. Now
we consider an $A_{6}^{8}(1)$\textbf{-}polytope given by $\{l_{i}\ 2\leq
i\leq8\}$. By proposition \ref{Prop-4cornerTo3+1Uncornered}, $\{l_{2}%
,l_{3},l_{4}\}\ $is a Fano block in $\{l_{i}\ 2\leq i\leq8\}$, and the
complement $\{l_{5},l_{6},l_{7},l_{8}\}$ gives a cornered $A_{3}^{8}%
(1)$\textbf{-}polytope. This gives the claim.

Here $D_{1}\cdot K_{S_{8}}=$ $D_{2}\cdot K_{S_{8}}=-4\ $and $D_{1}\cdot
D_{2}=16$. Since $D_{1}\ $and $D_{2}\ $are the centers of cornered $A_{3}%
^{8}(1)$\textbf{-}polytopes, these divisor classes are corresponded to
$l_{D_{1}}\ $and $l_{D_{2}}\ $by corollary \ref{CoroCornered-line}. Moreover,
the intersection $l_{D_{1}}\cdot l_{D_{2}}\ $is
\[
l_{D_{1}}\cdot l_{D_{2}}=\left(  \frac{1}{2}D_{1}+K_{S_{8}}\right)
\cdot\left(  \frac{1}{2}D_{2}+K_{S_{8}}\right)  =1\text{,}%
\]
and $l_{D_{1}}+l_{D_{2}}:=f\ $is a ruling. In fact, we obtain
\[
f=l_{D_{1}}+l_{D_{2}}=\left(  \frac{1}{2}D_{1}+K_{S_{8}}\right)  +\left(
\frac{1}{2}D_{2}+K_{S_{8}}\right)  =\frac{1}{2}D+2K_{S_{8}}\text{.}%
\]

For a fixed line $l_{8}$, each cornered $A_{3}^{8}(1)$\textbf{-}polytope
containing $l_{8}\ $gives a Fano block in $\{l_{i}\ 1\leq i\leq7\}$. Thus
there are seven cornered $A_{3}^{8}(1)$\textbf{-}polytopes containing
$l_{8}\ $with centers $D_{1}^{i}\ 1\leq i\leq7\ $by theorem \ref{Thm-Fano},
and more there are seven pairs of centers $D_{2}^{i}$ $1\leq i\leq7\ $of
cornered $A_{3}^{8}(1)$\textbf{-}polytopes such that $D_{1}^{i}+D_{2}^{i}%
=D\ $for $1\leq i\leq7$. Thus there are seven pairs of lines which produce the
ruling $f$. In fact, these are all the possible pairs in the ruling because
the ruling $f\ $is corresponding to a $7$-crosspolytope\ in $4_{21}$which has
seven pairs of bipolar lines with $1$-intersection.

$\ {\LARGE \blacksquare}$

\bigskip

As a direct application of the above theorem, we obtain the following
corollary which shows how to construct a $A_{7}^{8}(1)$-polytope from a
$A_{6}^{8}(1)$-polytope.

\begin{corollary}
Let $S:=\{l_{b}\ 1\leq b\leq7\}\ $give an $A_{6}^{8}(1)$-polytope\ where
$\{l_{i},l_{j},l_{k}\}\ $is a Fano block. If $l\ $is a line whose vertex
figure contains the cornered $A_{3}^{8}(1)$-polytope\ given by $S-\{l_{i}%
,l_{j},l_{k}\}$, then a $A_{3}^{8}(1)$-polytope\ given by $\{l,l_{i}%
,l_{j},l_{k}\}\ $is uncornered. Furthermore, there is a line $l_{8}$ whose
union with $S\ $(resp.$\{l_{i},l_{j},l_{k}\}$)\ gives an $A_{7}^{8}%
(1)$-polytope (resp. a cornered $A_{3}^{8}(1)$-polytope).
\end{corollary}

\textbf{Proof :\ }By Proposition \ref{Prop-7-Remain4-Uncornered}, the
$A_{3}^{8}(1)$-polytope with the center $D$\ given by $S-\{l_{i},l_{j}%
,l_{k}\}\ $is cornered, and the line $l\ $is given as $l=D/2+K_{S_{8}}$. By
lemma \ref{Lemma-3Skewgivesaline}, there is a line $l_{8}\ $such that
$\{l_{i},l_{j},l_{k},l_{8}\}\ $forms a cornered $A_{3}^{8}(1)$-polytope. Since
$\{l_{i},l_{j},l_{k}\}\ $is a Fano block in the $A_{3}^{8}(1)$-polytope,\ the
line $l_{8}\ $is not in $S$. Furthermore, by theorem \ref{Thm-uncorner-4to5},
$S\cup\{l_{8}\}\ $gives an $A_{7}^{8}(1)$-polytope. Since $l\cdot l_{a}=$
$\left(  D/2+K_{S_{8}}\right)  \cdot l_{a}=1\ $for $a=i,j,k$, $\{l,l_{i}%
,l_{j},l_{k}\}$ gives an $A_{3}^{8}(1)$-polytope, and this is uncornered by
proposition \ref{Prop-4cornerTo3+1Uncornered}. This gives the corollary.

${\LARGE \blacksquare}$

\bigskip

\subsection{Inscribed $2$-degree simplexes and $3$-degree simplexes}

\bigskip

As in the inscribed $1$-degree $m(\leq3)$-simplexes in $(r-4)_{21}$, the
centers of inscribed simplexes with higher (i.e. $>1$) degree$\ $are either
unique or corresponded to the lines for the cornered issues. Thus the
configurations of inscribed simplexes with higher degree in $(r-4)_{21}\ $are
naturally related to the $k$-Steiner systems along the monoidal transform.

\textbf{A. Inscribed }$2$\textbf{-degree }$1$\textbf{-simplex}

\bigskip The $A_{1}^{r}(2)$-polytopes exist when $r=7,8$.

\textbf{(1) Inscribed }$2$\textbf{-degree }$1$\textbf{-simplex\ in }$3_{21}$

Since a pair of lines $l_{1}\ $and $l_{2}\ $in\ $L_{7}\ $with $l_{1}\cdot
l_{2}=2\ $equivalently satisfy $l_{1}+l_{2}=$ $-K_{S_{7}}$, all $A_{1}^{7}%
(2)$-polytopes in $3_{21}$ share a center $-K_{S_{7}}$. Moreover, an
$A_{1}^{7}(2)$-polytope in $3_{21}\ $is determined by a line $l\ \ $in $L_{7}$
and $G(l)=-\left(  K_{S_{7}}+l\right)  $. Thus each $A_{1}^{7}(2)$%
-polytope$\ $corresponds to $2$-Steiner block on $S_{7}$ given by
$\{l,-\left(  K_{S_{7}}+l\right)  \}$.

\textbf{(2) Inscribed }$2$\textbf{-degree }$1$\textbf{-simplex in }$4_{21}$

The center of each $A_{1}^{8}(2)$-polytope$\ $in $Pic\ S_{8}$ is a divisor
$D:=l_{1}+l_{2}$ where $l_{1}\cdot l_{2}=2$. And as in the theorem
\ref{Thm-K8}, the divisor $D$ satisfies $D^{2}=\left(  l_{1}+l_{2}\right)
^{2}=2$ and $D\cdot K_{S_{8}}=\left(  l_{1}+l_{2}\right)  \cdot K_{S_{8}}=-2$.
The divisor $D$ on $S_{8}$ with $D^{2}=2$ and $D\cdot K_{S_{8}}=-2$
equivalently has $(D+K_{S_{8}})^{2}=-1$ and $(D+K_{S_{8}})\cdot K_{S_{8}}=-1$
which corresponds to a line in $Pic$ $S_{8}.$Thus the set of centers of
$A_{1}^{8}(2)$-polytopes$\ $in $4_{21}\ $is bijective to the $L_{8}$.

For a fixed center $D\ $of an $A_{1}^{8}(2)$-polytope$\ $in $4_{21}$, the
configuration of all the $A_{1}^{8}(2)$-polytopes with the common center
$D\ $is as follows.

Let $l_{3}\ $and $l_{4}$ be two lines of an $A_{1}^{8}(2)$-polytope$\ $with
the center $D=l_{3}+l_{4}$ above, we have a line $l:=D+K_{S_{8}}$ which also
satisfies $\left(  D+K_{S_{8}}\right)  \cdot l_{i}=\left(  l_{3}%
+l_{4}+K_{S_{8}}\right)  \cdot l_{i}=0\ $for $i=3,4$. Therefore, $l_{3}\ $and
$l_{4}\ $are in $N_{0}(l,S_{8})$, namely, the vertex figure of $l$. The map
$\pi_{l}^{8}:S_{8}\rightarrow S_{7}\ $induced by blowing down $l$ sends
$l_{3}\ $and $l_{4}\ $to lines in $S_{7}\ $with intersection $2$. Thus
configuration of $A_{1}^{8}(2)$-polytopes$\ $with the fixed center $D\ $is
same with the configuration of $A_{1}^{7}(2)$-polytopes in $3_{21}\ $which is
$2$-Steiner system $\mathcal{S}_{A}(2,S_{7})$ . Note by lemma
\ref{Lemma-3Skewgivesaline} the $G_{l}$, Gieser transform of $l$, transfer
$l_{1}\ $and $l_{2}\ $to each other.

\bigskip

The configuration of inscribe $2$-degree $1$-simplexes are summarized as follows.

\begin{theorem}
The $2$-Steiner system$\ \mathcal{S}_{A}(2,S_{7})\ $in $S_{7}\ $determines the
configuration of $A_{1}^{8}(2)$-polytopes in $4_{21}$\ with a fixed center and
the configuration of $A_{1}^{7}(2)$-polytopes$\ $in $3_{21}$.
\end{theorem}

\bigskip

\textbf{B. Inscribed }$2$\textbf{-degree }$2$\textbf{-simplex in }$4_{21}$

The $A_{2}^{8}(2)$-polytopes exist when $r=8$. The center of each $A_{2}%
^{8}(2)$-polytopes$\ $in $Pic\ S_{8}$ is represented by a divisor
$D:=l_{1}+l_{2}+l_{3}$ where $l_{1},l_{2}$ and $l_{3}$ lines with $l_{i}\cdot
l_{j}=2$. And the divisor $D$ satisfies $D^{2}=\left(  l_{1}+l_{2}%
+l_{3}\right)  ^{2}=9$ and $D\cdot K_{S_{8}}=\left(  l_{1}+l_{2}+l_{3}\right)
\cdot K_{S_{8}}=-3$. By Hodge index theorem, $-3K_{S_{8}}$ is the only divisor
satisfying equations. Thus, all the $A_{2}^{8}(2)$-polytopes in $4_{21}$ share
one center $-3K_{S_{8}}$.

Furthermore, each $A_{2}^{8}(2)$-polytope in $4_{21}$ corresponds to a
$3$-Steiner block of $\mathcal{S}_{B}(3,S_{8})\ $in the section
\ref{Sec-Steiner}. Therefore, the configuration of the $A_{2}^{8}%
(2)$-polytopes in $4_{21}\ $is exactly the $3$-Steiner system $\mathcal{S}%
_{B}(3,S_{8})$.

\bigskip

\textbf{C. Inscribed }$3$\textbf{-degree }$1$\textbf{-simplexes in }$4_{21}$

The $A_{1}^{8}(3)$-polytopes exist when $r=8$. Since a pair of lines $l_{1}%
\ $and $l_{2}\ $in\ $L_{8}\ $with $l_{1}\cdot l_{2}=3\ $equivalently satisfy
$l_{1}+l_{2}=$ $-2K_{S_{8}}$, all $A_{1}^{8}(3)$-polytopes in $4_{21}$ share a
center $-2K_{S_{8}}$. Moreover, an $A_{1}^{8}(3)$-polytope in $4_{21}\ $is
determined by a line $l\ \ $in $L_{8}$ and $B(l)=-\left(  2K_{S_{8}}+l\right)
$. Thus each$\ A_{1}^{8}(3)$-polytope$\ $corresponds to $2$-Steiner block in
$\mathcal{S}_{A}(2,S_{8})$ given by $\{l,-\left(  2K_{S_{8}}+l\right)  \}$ in
$S_{8}$.

\bigskip

In summary of this section, we have following tables.%

\[
\underset{\text{Number of the Centers of\ }A_{a}^{r}(1)\text{-simplex in
}(r-4)_{21}\ r\leq7\text{ }}{%
\begin{tabular}
[c]{|l|l|l|l|l|l|}\hline
& $-1_{21}$ & $0_{21}$ & $1_{21}$ & $2_{21}$ & $3_{21}$\\\hline\hline
$A_{1}(1)$ & $3$ & $5$ & $10$ & $27$ & $126$\\\hline
$A_{2}(1)$ &  &  &  & $1$ & $56$\\\hline
$A_{3}(1)$ &  &  &  &  & $1$\\\hline
\end{tabular}
\ \ \ }%
\]%
\[
\underset{\text{Number of the Centers of }A_{k}^{r}(m)\ m>1}{\
\begin{tabular}
[c]{|l|l|l|}\hline
& $3_{21}$ & $4_{21}$\\\hline\hline
$A_{1}(2)$ & $1$ & $240$\\\hline
$A_{2}(2)$ &  & $1$\\\hline
$A_{1}(3)$ &  & $1$\\\hline
\end{tabular}
\ \ \ }%
\]

\[
\underset{\text{Number of the Centers of }A_{k}^{8}(1)\text{ in }4_{21}}{%
\begin{tabular}
[c]{|l|l|l|l|l|l|l|}\hline
$A_{1}(1)$ & $A_{2}(1)$ & $A_{3}(1)$ & $A_{4}(1)$ & $A_{5}(1)$ & $A_{6}(1)$ &
$A_{7}(1)$\\\hline\hline
$2160$ & $6720$ & $17520$ & $30240$ & $60480$ & $207360$ & $2160$\\\hline
\end{tabular}
\ }%
\]

\subsection{Inscribed crosspolytopes in Gosset polytopes}

\bigskip

Recall that each ruling in $Pic\ S_{r}\ $corresponds to a center of an
$(r-1)$-crosspolytope in $(r-4)_{21}$. The configuration of the lines in the
$(r-1)$-crosspolytope\ is given by $(r-1)$-pair of lines with $1$-intersection
between the two lines in each pair and $0$-intersection for other cases (see
subsection \ref{Subsec-Gosset-Picard}). Furthermore, when we choose a line
from each of $(r-1)$-pairs, the $(r-1)\ $chosen lines always consists of an
$(r-1)$-simplex in the $(r-1)$-crosspolytope, and two lines in a pair and
another line form an isosceles right triangle with lengths $\sqrt{2}$,
$\sqrt{2}\ $and $2$. Here the common length of edges in the $(r-1)$%
-crosspolytope\ is $\sqrt{2}$. In this subsection, we consider inscribed
$m$-crosspolytopes in $(r-4)_{21}\ $for\ $m\geq2,\ $which have the edges with
common length $>\sqrt{2}$.

Recall that for lines $l_{1}\ $and $l_{2}\ $in $L_{r}$, the length of the
segment joining $l_{1}\ $and $l_{2}\ $is given by$\sqrt{-(l_{1}-l_{2})^{2}%
}=\sqrt{2+2l_{1}\cdot l_{2}}$. If each edges in a crosspolytope has length
$2$, then the distance between two antipodal lines in the crosspolytope is
$2\sqrt{2},\ $and their intersection is $3$. Two lines with intersection
$3\ $exist only for $S_{8}$, and since there are no two lines in a del Pezzo
surface whose intersection is bigger than $3$, this is only possible case for
the inscribed crosspolytopes in $(r-4)_{21}\ $for\ $m\geq2$.

Now, we only consider lines in $L_{8}\ $and a Gosset polytope $4_{21}$.

Just like the center of an $(r-1)$-crosspolytope in $(r-4)_{21}$, the center
$D\ $of an inscribed crosspolytope is given by the sum of two antipodal lines
in it. Since these lines have $3$-intersection, the divisor class
$D\ $satisfies $D^{2}=-8$,\ $D\cdot K_{S_{8}}=-2$. Because $-2K_{S_{8}}$ is
only $D\ $with these conditions, all the inscribed crosspolytopes in
$(r-4)_{21}\ $share a common center $-2K_{S_{8}}$.

If there is an inscribed $m$-crosspolytope in $(r-4)_{21}\ $for\ $m\geq2$, the
facets in it are $(m-1)$-simplexes with length $2$, namely, inscribed
$1$-degree $(m-1)$-simplexes in $4_{21}$. Conversely, for each inscribed
$1$-degree $(m-1)$-simplex in $4_{21}$, we can construct an inscribed
$m$-crosspolytope as follows.

\begin{theorem}
\label{Thm-crosspolytope}Let $l_{i}\ 1\leq i\leq m\ $be lines in $L_{8}%
\ $which form an $A_{m-1}^{8}(1)$-polytope in $4_{21}$. The set $\{l_{i}%
,B(l_{i})\ 1\leq i\leq m\}\ $gives an inscribed $m$-crosspolytope in
$4_{21}\ $where one of its facet is the $A_{m-1}^{8}(1)$-polytope.
Furthermore, there are inscribed $m$-crosspolytope in $4_{21}\ $for $m\leq8.$
\end{theorem}

\bigskip

The inscribed $m$-crosspolytopes in $4_{21}\ $are involved in many interesting
issues along the configuration of lines. However, we consider the
$4$-crosspolytopes in $4_{21}\ $in the next section, and the other issues will
be discussed in another article.

\bigskip

\section{Hypercubes in Gosset polytopes}

In this section, we consider polytopes given by both the vertices and the
edges of $(r-4)_{21}\ $but inscribed in it. In particular we study $m$-cubes
in $(r-4)_{21}$. In fact, the inscribed $3$-cubes in $3_{21}$ was studied by
Coxeter \cite{Coxeter2}.

An $m$\textit{-cube} $\gamma_{m}\ $(also called as an $m$-dimensional
hypercube or a measure polytope) is a parallelogram given by perpendicular
$m$-vectors of a constant length, which represent $m$-edges that meet at one
vertex. For example, $\gamma_{2}\ $is a square and $\gamma_{3}\ $is a
cube.\ An $m$-cube can be constructed by translating an $(m-1)$-cube along a
segment in the perpendicular direction. An $m$-cube\ contains $2^{m}%
$-vertices, and the vertex figure of each vertex in it is a regular
$(m-1)$-simplex $\alpha_{m-1}$, see \cite{Coxeter} for detail.

On behalf of figuring out possible $m$-cube in a Gosset polytope $(r-4)_{21}$,
we observe that the length of a diagonal in the $m$-cube can not be bigger
than the longest distance between vertices in $(r-4)_{21}$. Since each edge in
the cube has length $\sqrt{2}$, a diagonal in an $m$-cube is $\sqrt{2m}$, and
the longest distance between vertices in $(r-4)_{21}\ $is $2\ $for $3\leq
r\leq6$, $\sqrt{6}\ $for $r=7\ $and $\sqrt{8}\ $for $r=8\ $because
$(l_{1}-l_{2})^{2}=-2-2l_{1}\cdot l_{2}$ and the maximal possible
intersections in $L_{r}\ $is $1\ $for $3\leq r\leq6$, $2\ $for $r=7\ $and
$3\ $for $r=8$. Thus the dimension of biggest $m$-cubes in $(r-4)_{21}$ is
$2\ $for $3\leq r\leq6$, $3\ $for $r=7\ $and $4\ $for $r=8$.

The $2$-cubes are given by two lines with $1$-intersection. Thus $2$-cubes in
$(r-1)_{21}\ $correspond to rulings in $Pic\ S_{r}$ which is treated in
\cite{Lee}. In this article, we only study $3$-cubes in $3_{21}\ $and
$4$-cubes in $4_{21}\ $as an application of the above study on inscribed
simplexes and crosspolytopes, and the further study on $m$-cubes in
$(r-4)_{21}\ $will be discussed in another article.

For an $m$-cube, the center is given by the sum of two vertices on a diagonal
in the cube. Thus the center of an $m$-cube equals the center of an inscribed
$(m-1)$-degree $1$-simplex. For each vertex $l$ in $m$-cube, the vertices in
its vertex figure intersects by $1\ $to each other because any two lines in
the vertex figure of $l\ $are two diagonal vertices of an $2$-cube containing
a vertex $l\ $. \ Thus the vertex figure of $l\ $in the $m$-cube is an
inscribed $1$-degree $(m-1)$-simplex.

\subsection{3-cubes in $3_{21}\ $and $4_{21}$}

\bigskip

\textbf{A. }$3$\textbf{-cubes in }$3_{21}$

The centers of $3$-cubes in $3_{21}\ $are corresponded to the centers of
$A_{1}^{7}(2)$-polytopes in $3_{21}$. Thus those are given by the divisor
classes $D\ $with\ $D^{2}=$ $2\ $and $D\cdot K_{S_{7}}=-2$. Since such $D\ $in
$Pic\ S_{7}\ $is only $-K_{S_{7}}$. All the $3$-cubes in $3_{21}\ $share a
common center. And the vertex figures of lines in the $3$-cubes are $A_{2}%
^{7}(1)$-polytopes in $3_{21}$. In fact, an $A_{2}^{7}(1)$-polytope determines
a $3$-cube as follows.

We consider an $A_{2}^{7}(1)$-polytope in $3_{21}\ $given by $\{l_{1}%
,l_{2},l_{3}\}\ $with the center $D_{1}=l_{1}+l_{2}+l_{3}$. As in the
configuration of $A_{2}^{7}(1)$-polytopes in $3_{21}$, there is a line
$l_{D_{1}}^{S_{7}}=D_{1}+K_{S_{7}}\ $in $L_{7}\ $where $\{l_{1},l_{2}%
,l_{3}\}\ $is in the vertex figure of $l_{D_{1}}^{S_{7}}$. By applying Gieser
transform $G$, we get another $A_{2}^{7}(1)$-polytope in $3_{21}\ $given by
$\{G(l_{1}),G(l_{2}),G(l_{3})\}$ and a line $G(l_{D_{1}})\ $where its vertex
figure contains $\{G(l_{1}),G(l_{2}),G(l_{3})\}$. Now, a set of lines
\[
\{l_{1},l_{2},l_{3},l_{D_{1}},G(l_{1}),G(l_{2}),G(l_{3}),G(l_{D_{1}})\}\
\]
gives a $3$-cube in $3_{21}\ $because $l_{i}\cdot G(l_{j})=l_{i}\cdot\left(
-K_{S_{7}}-l_{j}\right)  =0\ $for $i\not =j$. This also implies that the set
of lines in a $3$-cube in $3_{21}\ $can be written this way. In fact, by the
remark in the inscribed $1$-degree $2$-simplexes in $3_{21}$,\ we know
$l_{i}+l_{j}=l_{D}^{S_{7}}+G(l_{k})$,\ $i,j,k\ $distinct and $\{l_{i}%
,l_{j},l_{D}^{S_{7}},G(l_{k})\}$ forms a $2$-cube.

\bigskip

\textbf{B. }$3$\textbf{-cubes in }$4_{21}$

The centers of $3$-cubes in $4_{21}\ $are corresponded to the centers of
$A_{1}^{8}(2)$-polytopes in $4_{21}$. Thus those are given by the divisor
classes $D\ $with\ $D^{2}=$ $2\ $and $D\cdot K_{S_{8}}=-2$. As in the
configuration of $A_{1}^{8}(2)$-polytopes, such $D\ $in $Pic\ S_{8}\ $is
corresponded to a line in $L_{8}\ $by $l_{D}=D+K_{S_{8}}$ where all the
$A_{1}^{8}(2)$-polytopes in $4_{21}$ with the center $D\ $are contained in the
vertex figure of $l_{D}$. Thus the $3$-cubes in $4_{21}\ $with the center
$D\ $are in the vertex figure of $l_{D}$. Furthermore, these $3$-cubes in
$4_{21}\ $with the center $D\ $is bijectively mapped to the $3$-cubes in
$3_{21}\ $by the blow down map $\pi_{l_{D}}:S_{8}\rightarrow S_{7}$. Just like
$3$-cube in $3_{21}$, a $3$-cube with fixed center $D$ in $4_{21}\ $is
determined by a vertex and its vertex figure in the $3$-cube by the Gieser
transform for $l_{D}\ $in $4_{21}$. In other words, if $\{l_{1},l_{2}%
,l_{3}\}\ $is an $A_{2}^{8}(1)$-polytope in the vertex figure in a $3$-cube
with center $D$, the $3$-cube is given by
\[
\{l_{1},l_{2},l_{3},l^{\prime},G_{l_{D}}(l_{1}),G_{l_{D}}(l_{2}),G_{l_{D}%
}(l_{3}),G_{l_{D}}(l^{\prime})\}\
\]
where $l^{\prime}:=l_{1}+l_{2}+l_{3}+K_{S_{8}}-l_{D}$\ by lemma
\ref{Lemma-3Skewgivesaline}.

\bigskip

\subsection{4-cubes in $4_{21}$}

The centers of $4$-cubes in $4_{21}\ $are corresponded to the centers of
$A_{1}^{8}(3)$-polytopes in $4_{21}$. Thus those are given by the divisor
classes $D\ $with\ $D^{2}=$ $4\ $and $D\cdot K_{S_{8}}=-2$. Since such $D\ $in
$Pic\ S_{8}\ $is only $-2K_{S_{8}}$. All the $4$-cubes in $4_{21}\ $share a
common center. And the vertex figures of lines in the $4$-cubes are cornered
$A_{3}^{8}(1)$-polytopes in $4_{21}$. In the following, a cornered $A_{3}%
^{8}(1)$-polytope determines a $4$-cube in $4_{21}$. We also consider an
uncornered $A_{3}^{8}(1)$-polytope along the $4$-cubes in $4_{21}$.

\bigskip

\textbf{A. }$4$\textbf{-cubes}$\ $\textbf{and inscribed cornered }%
$1$\textbf{-degree }$3$\textbf{-simplexes in }$4_{21}$

Recall a $4$-cube has $2^{4}\ $vertices. Since we begin with a $3$-simplex, we
need to identify $12\ $more lines to construct a $4$-cube.

We consider a cornered $A_{3}^{8}(1)$-polytope in $4_{21}\ $given by the set
of lines $\{l_{1},l_{2},l_{3},l_{4}\}\ $in $L_{8}$, and $l\ $is the line in
$L_{8}\ $whose vertex figure contains $\{l_{1},l_{2},l_{3},l_{4}\}$. By
applying Bertini transform $B$, another set of lines $\{B(l_{1}),B(l_{2}%
),B(l_{3}),B(l_{4})\}\ $in $L_{8}\ $produces a cornered $A_{3}^{8}%
(1)$-polytope which is in the vertex figure of $B(l)$.\ Now, we need to
identify six$\ $more lines to construct a $4$-cube. The idea is that we choose
proper four lines from
\[
\{l_{1},l_{2},l_{3},l_{4},B(l_{1}),B(l_{2}),B(l_{3}),B(l_{4})\}
\]
$\ $so that four chosen lines form a cornered $A_{3}^{8}(1)$-polytope. This
$3$-simplex determine a line whose vertex figure contains the $3$-simplex.
According to the intersections between lines, $\{l_{i},l_{j},l_{k},B(l_{m}%
)\}$, $\{B(l_{i}),B(l_{j}),B(l_{k}),l_{m}\}\ $and $\{l_{i},l_{j}%
,B(l_{k}),B(l_{m})\}\ $where $\{i,j,k,l\}=\{1,2,3,4\}\ $are the only possible
cases. But by proposition \ref{Prop-4cornerTo3+1Uncornered}, first two cases
produce uncornered $A_{3}^{8}(1)$-polytopes. For the last case, there are
$6\ $possible candidates, and all of these produce cornered $A_{3}^{8}%
(1)$-polytopes by following lemma.

\begin{lemma}
If the set of lines $\{l_{1},l_{2},l_{3},l_{4}\}\ $in $L_{8}\ $produces a
cornered $A_{3}^{8}(1)$-polytope in $4_{21}$, each set of lines $\{l_{i}%
,l_{j},B(l_{k}),B(l_{m})\}\ $for $\{i,j,k,l\}=\{1,2,3,4\}\ $gives a cornered
$A_{3}^{8}(1)$-polytope in $4_{21}$.
\end{lemma}

\textbf{Proof :} By proposition \ref{Prop-root-3simplex}, it suffices to show
$l_{i}+l_{j}+B(l_{k})+B(l_{m})+4K_{S_{8}}\ $is $2d^{\prime}\ $for a
root$\ d^{\prime}\ $in $Pic\ S_{8}$.

Since $\{l_{1},l_{2},l_{3},l_{4}\}\ $gives a cornered $A_{3}^{8}(1)$-polytope,
there is a root $d\ \ $such that $2d=l_{1}+l_{2}+l_{3}+l_{4}+4K_{S_{8}}$. For
$\{l_{i},l_{j},B(l_{k}),B(l_{m})\}$, we have
\begin{align*}
l_{i}+l_{j}+B(l_{k})+B(l_{m})+4K_{S_{8}}  &  =l_{i}+l_{j}-l_{k}-l_{m}\\
&  =2(l_{i}+l_{j}-d+2K_{S_{8}}).
\end{align*}
Moreover, $l_{i}+l_{j}-d-2K_{S_{8}}\ $is a root because%
\[
\left(  l_{i}+l_{j}-d+2K_{S_{8}}\right)  \cdot K_{S_{8}}=0\ \text{and }\left(
l_{i}+l_{j}-d+2K_{S_{8}}\right)  ^{2}=-2\text{.}%
\]
Thus, we have the lemma.$\ \ \ \ \ \ \ \ \ \ \ \ \ {\LARGE \blacksquare}$

\bigskip

By this lemma, we have six$\ $more lines in addition to
\[
\{l_{1},l_{2},l_{3},l_{4},l,B(l_{1}),B(l_{2}),B(l_{3}),B(l_{4}),B(l)\}
\]
which complete $16\ $lines for a $4$-cube in $4_{21}$. In fact, we have an
inscribed $1$-degree $4$-crosspolytope given by
\[
\{l_{1},l_{2},l_{3},l_{4},B(l_{1}),B(l_{2}),B(l_{3}),B(l_{4})\}
\]
$\ $and a $4$-cube containing it. In summary, we have the following theorem.

\bigskip

\begin{theorem}
Let the set of lines $\{l_{1},l_{2},l_{3},l_{4}\}\ $in $L_{8}\ $produce a
cornered $A_{3}^{8}(1)$-polytope in $4_{21}\ $and $l\ \ $be the line whose
vertex figure containing the $3$-simplex. A set of lines given by
\[
\left\{
\begin{array}
[c]{c}%
l,B(l),\ l_{i},B(l_{i})\ \text{for }1\leq i\leq4,\ \\
\frac{1}{2}(l_{i}+l_{j}+B(l_{k})+B(l_{m}))+K_{S_{8}}\ \text{for\ }%
\{i,j,k,l\}=\{1,2,3,4\}
\end{array}
\right\}
\]
makes a $4$-cube in $4_{21}$. Furthermore, an inscribed $1$-degree
$4$-crosspolytope given by $\{l_{i},B(l_{i})\ (1\leq i\leq4)\}\ $is inscribed
in the $4$-cube.\ 
\end{theorem}

\bigskip

\textbf{B.\ Inscribed }$1$\textbf{-degree }$4$\textbf{-crosspolytopes and
inscribed uncornered }$1$\textbf{-degree }$3$-\textbf{simplexes in }$4_{21}$

When the set of lines $\{l_{1},l_{2},l_{3},l_{4}\}\ $in $L_{8}\ $produce an
uncornered $A_{3}^{8}(1)$-polytope in $4_{21}$, we can also construct an
inscribed $1$-degree $4$-crosspolytope in $4_{21}$.\ But there is no $4$-cube
containing the $4$-crosspolytope.

From lemma \ref{Lemma-3Skewgivesaline}, each $A_{3}^{8}(1)$-polytope given by
$\{l_{i},l_{j},l_{k}\}\ $in the $3$-simplex\ gives a skew $2$-lines
$l_{a}^{ijk}+l_{b}^{ijk}\ $where $l_{i}+l_{j}+l_{k}+K_{S_{8}}=$\ $l_{a}%
^{ijk}+l_{b}^{ijk}\ $for $\{i,j,k,l\}=\{1,2,3,4\}$. Furthermore, $l_{a}%
^{ijk}\cdot l_{b}^{lmn}=l_{a}^{ijk}\cdot l_{a}^{lmn}=l_{b}^{ijk}\cdot
l_{b}^{lmn}=0\ $for $\{i,j,k\}\not =\{l,m,n\}\ $since above four skew
$2$-lines are skew edges. Now we consider a set of eight lines given by%
\[
C:=\{G_{l_{a}^{ijk}}(l_{b}^{ijk}),G_{l_{b}^{ijk}}(l_{a}^{ijk})\text{ where
}ijk\ \in\{123,124,134,234\}\}\text{.}%
\]
By lemmar \ref{Lemma-3Skewgivesaline}, $G_{l_{a}^{ijk}}(l_{b}^{ijk})\cdot
G_{l_{a}^{lmn}}(l_{b}^{lmn})=G_{l_{a}^{ijk}}(l_{b}^{ijk})\cdot G_{l_{b}^{lmn}%
}(l_{a}^{lmn})=1\ $for $\{i,j,k\}\not =\{l,m,n\}\ $and $G_{l_{a}^{ijk}}%
(l_{b}^{ijk})\cdot G_{l_{b}^{ijk}}(l_{a}^{ijk})=3$. Thus, the set $C\ $is the
vertex set of an inscribed $1$-degree $4$-crosspolytope in $4_{21}$. From the
following theorem, there is no $4$-cube containing the inscribed $4$-crosspolytope.

\begin{theorem}
Let the set of lines $\{l_{1},l_{2},l_{3},l_{4}\}\ $in $L_{8}\ $produce an
uncornered $A_{3}^{8}(1)$-polytope in $4_{21}$. An inscribed $1$-degree
$4$-crosspolytope in $4_{21}\ $determined as above is not inscribed in any
$4$-cube in $4_{21}$.
\end{theorem}

\textbf{Proof:} We show that each $A_{3}^{8}(1)$-polytopes in the
$4$-crosspolytope is uncornered so that none of the facet in the
$4$-crosspolytope\ can be a vertex figure of an $4$-cube. We consider an
$A_{3}^{8}(1)$-polytope given by
\[
\ S:=\{G_{l_{a}^{ijk}}(l_{b}^{ijk})\ \ \text{where }ijk\ \in
\{123,124,134,234\}\}
\]
$\ $and assume it is cornered. Thus there is a root $d\ $such that
\[
2d=G_{l_{a}^{123}}(l_{b}^{123})+G_{l_{a}^{124}}(l_{b}^{124})+G_{l_{a}^{134}%
}(l_{b}^{134})+G_{l_{a}^{234}}(l_{b}^{234})+4K_{S_{8}}.
\]
Furthermore, this gives
\[
2d=(l_{a}^{123}+l_{a}^{124}+l_{a}^{134}+l_{a}^{234})-(l_{b}^{123}+l_{b}%
^{124}+l_{b}^{134}+l_{b}^{234}).
\]
On the other hand, since $\{l_{1},l_{2},l_{3},l_{4}\}\ $forms an uncornered
$A_{3}^{8}(1)$-polytope, by corollary \ref{Coro-Skew8-cornered} the center of
this $3$-simplex satisfies
\[
3(l_{1}+l_{2}+l_{3}+l_{4})+4K_{S_{8}}=\sum_{i=a,b}(l_{i}^{123}+l_{i}%
^{124}+l_{i}^{134}+l_{i}^{234})\text{.}%
\]
By combining above two equalities, we have
\[
3(l_{1}+l_{2}+l_{3}+l_{4})=2\{(l_{a}^{123}+l_{a}^{124}+l_{a}^{134}+l_{a}%
^{234})-d-2K_{S_{8}}\}\text{,}%
\]
which implies $\frac{1}{2}(l_{1}+l_{2}+l_{3}+l_{4})\ $in $Pic\ S_{8}$\ and it
is a root in $S_{8}$. By proposition \ref{Prop-root-3simplex}, the $A_{3}%
^{8}(1)$-polytope given by\ $\{l_{1},l_{2},l_{3},l_{4}\}\ $must be cornered,
and this is a contradiction. Thus the $A_{3}^{8}(1)$-polytope given by the set
$S\ $is uncornered. Similarly, each $A_{3}^{8}(1)$-polytopes in the
$4$-crosspolytope\ is uncornered.

This gives the theorem.$\ \ {\LARGE \blacksquare}$

\bigskip

\bigskip

\textit{Acknowledgments: Author is always grateful to his mentor Naichung
Conan Leung for everything he gives to the author. Author thanks to Adrian
Clingher for useful discussions.}

{\scriptsize Addresses: }

{\scriptsize Jae-Hyouk Lee (lee@math.umsl.edu)}

{\scriptsize Department of Mathematics and Computer Science, University of
Missouri-St. Louis, U.S.A.}

\end{document}